\documentclass[11pt]{article}
\usepackage{latexsym}
\usepackage{amssymb} 

\input {xy}
\usepackage{amsmath,amsthm,amsfonts,amssymb}
\usepackage[matrix,arrow]{xy}

\usepackage{graphicx}
\newtheorem{teo}{Theorem}[section]
\newtheorem{coro}[teo]{Corollary}
\newtheorem{lema}[teo]{Lemma}
\newtheorem{prop}[teo]{Proposition}
\newtheorem{defi}[teo]{Definition}
\newtheorem{obs}[teo]{Remark}
\newtheorem{ejem}[teo]{Example}

\usepackage{fullpage}

\def\id{{\mathrm{id}}}
\def\cl#1{{\langle #1 \rangle}}

\def\qed{\begin{flushright} \QED \end{flushright}}

\newcommand\ZZ{{\mathbb{Z}}}
\newcommand\NN{{\mathbb{N}}}

\def\A{{\mathcal A}}
\def\B{{\mathcal B}}
\def\C{{\mathcal C}}
\def\D{{\mathcal D}}
\def\E{{\mathcal E}}
\def\F{{\mathcal F}}
\def\I{{\mathcal I}}
\def\J{{\mathcal J}}
\def\Z{{\mathcal Z}}

\def\place{{-}}

\def\qed{\hfill \mbox{$\square$}\bigskip}

\begin{document}

\sf  \title{Hochschild-Mitchell cohomology and Galois extensions}
\author{Estanislao Herscovich and Andrea Solotar 
\thanks{This work has been supported by the projects PICT 08280 (ANPCyT), UBACYTX169 and PIP-CONICET 02265. 
The first author is a CONICET fellow. 
The second author is a research member of CONICET (Argentina) and a Regular 
Associate of ICTP Associate Scheme.}}

\date{}  

\maketitle  

\begin{abstract}
We define $H$-Galois extensions for $k$-linear categories and prove the existence of a Grothendieck spectral sequence for Hochschild-Mitchell 
cohomology related to this situation. 
This spectral sequence is multiplicative and for a group algebra decomposes as a direct sum indexed by conjugacy classes of the group. 
We also compute some Hochschild-Mitchell cohomology groups of categories with infinite associated quivers. 
\end{abstract} 



\section{Introduction}

Hochschild cohomology of finite dimensional algebras over a field $k$ has been studied using the fact that a finite dimensional $k$-algebra is Morita 
equivalent to an algebra $kQ/I$ where $Q$ is a finite quiver and $I$ is an admissible ideal. 
An important object related to this kind of algebras is its universal Galois covering, which has played a central role since its introduction by \cite{B-G1}, 
\cite{D-S1}, and \cite{Ga1}. 
It is useful, for example, for the study of the representation type of the algebra or its fundamental group. 

One of the difficulties arising in this context when one needs to use cohomological methods is that this universal Galois covering is not in general a 
$k$-algebra but a $k$-linear category. 
This problem may be solved, as Cibils and Redondo did in \cite{C-R1}, using Hochschild-Mitchell cohomology, instead of Hochschild cohomology 
which is only defined for $k$-algebras. 
Another possible approach may be to use Hochschild cohomology for non-unital algebras, but we think that it is not the best one since the category of 
bimodules over a non-unital algebra is just an epimorphic image of the category of modules over the $k$-linear category associated to this algebra. 

The Hochschild-Mitchell cohomology of a $k$-linear category has been defined by Mitchell in \cite{Mit1} and is related with cohomology theories studied 
by Keller and McCarthy (cf. \cite{Kel1}, \cite{Mc1}). 
A wider context is studied in \cite{P-R1}. 

Cibils and Redondo proved in \cite{C-R1} that given a universal Galois covering $\C$ of a $k$-linear category $\B$ with Galois group $G$ and a 
$\B$-bimodule $M$ there is a spectral sequence $H^{\bullet}(G , H^{\bullet}(\C , LM))$ converging to $H^{\bullet}(\B , M)$ where $LM$ is the 
induced $\C$-bimodule. 
This result is a particular case of the spectral sequence constructed in \cite{P-R1}. 

This article has two main purposes: 

Firstly, we prove that the spectral sequence of \cite{C-R1} is multiplicative and that it can be decomposed into a direct sum of subsequences 
indexed by the set of conjugacy classes of $G$. 
This result is proved using Grothendieck spectral sequences and using 
Morita equivalences (see \cite{C-S1} for a complete description of 
Morita equivalence of $k$-linear categories). 

In particular we provide a proof of Morita invariance of Hochschild-Mitchell 
(co)homology. 

In fact, we extend the situation of \cite{C-R1} to $H$-module categories, for any Hopf algebra $H$ acting on $\C$. 
All these data fit into the situation studied by Sarah Witherspoon, so it is possible to describe the mulplicative structure (cf. \cite{Wi1}). 

One of the advantages of our approach is that one can also take into account twisted smash products (by a $2$-Hochschild-Mitchell cocycle) obtaining 
in this way coverings of the algebra which are not coverings with the usual 
definition.

Secondly, we compute the Hochschild-Mitchell cohomology of some linear categories of type $kQ/I$ where $Q$ is an infinite quiver and $I$ is an 
admissible ideal. 
Examples of such categories are: trees, radical square zero categories and other categories arising from posets. 

This article is organized as follows: 

In the second section we recall some definitions and properties concerning $k$-linear categories and Hochschild-Mitchell cohomology. 

In the third section we construct the spectral sequence in a wider context, namely for $H$-module categories ($H$ a Hopf algebra). 
We also prove that it is multiplicative and we prove the decomposition indexed by conjugacy classes of the Galois group $G$. 

In the fourth section the decomposition of (co)homology in terms of
the set of conjugacy classes of $G$ is proven.

Finally, in the fifth section we compute the Hochschild-Mitchell cohomology 
groups of categories associated to infinite quivers. 

All unadorned tensor products will be over the field $k$, i.e, $\otimes = \otimes_{k}$.

\section{Preliminaries}
\label {defyresbas}

In this section we recall some basic definitions and facts about Hochschild-Mitchell (co)homology. 
For further references, see  \cite{Mit2}, \cite{C-M1} and \cite{C-R1}.    

Let us consider a field $k$ and a small category $\C$. 
 
$\C$ is a \textbf{$k$-linear category} if the set of morphisms between two arbitrary objects of $\C$ is a $k$-vector space and 
composition of morphisms is $k$-bilinear. 
From now on, $\C$ will be a $k$-linear category with set of objects  $\C_{0}$ and given objects $x, y$ we shall denote ${}_{y}\C_{x}$ the 
$k$-vector space of morphisms from $x$ to $y$ in $\C$. 
Given $x,y,z$ in $\C_{0}$, the composition is a $k$-linear map 
\[     \circ_{z,y,x} : {}_{z}\C_{y} \otimes {}_{y}\C_{x} \to {}_{z}\C_{x}.     \]
We shall denote ${}_{z}f_{y}.{}_{y}g_{x}$, ${}_{z}f_{y}\cdot{}_{y}g_{x}$, ${}_{z}f_{y} \circ {}_{y}g_{x}$, or $f g$ if subscripts are clear, the image of 
${}_{z}f_{y} \otimes {}_{y}g_{x}$ under this map.
By \textbf{cardinality of the category $\C$} we will always mean the cardinality of $\C_{0}$. 
A $k$-linear category is called of \textbf {locally finite dimension} if ${}_{y}\C_{x}$ is of finite dimension for any $x,y \in \C_{0}$.

The simplest example of $k$-linear category is to look at a $k$-algebra $A$ as a category with only one object and the set of morphisms equal to $A$.  

\begin{obs}
We use the term locally finite dimension instead of locally finite, used in \cite{C-R1}, since this last term has been used by Takeuchi with a 
different meaning (see \cite{Ta1}).
\end{obs}

\begin{defi}
Given two $k$-linear categories $\C$ and $\D$ \textbf{the (external) tensor product category}, which we denote $\C \boxtimes_{k} \D$, is the category with 
set of objects $\C_{0} \times \D_{0}$ and given $c, c' \in \C_{0}$ and $d,d' \in \D_{0}$ 
\[     {}_{(c',d)}(\C \boxtimes_{k} \D)_{(c,d)} = {}_{c}\C_{c} \otimes {}_{d'}\D_{d}.       \]

The functor $\C \boxtimes_{k} \place$ is the left adjoint functor to $\mathrm{Func}(\C, \place)$ (see \cite{Mit2}, section 2, p. 13; or \cite{Hers1}, p. 28, for a 
detailed proof). 
We will omit the subindex $k$ in the external tensor product. 
We will call the category $\C \boxtimes \C^{op}$ the \textbf{enveloping category of $\C$} and denote it $\C^{e}$. 
\end{defi}

\begin{defi}
A \textbf{left $\C$-module} $M$ is a covariant $k$-linear functor from the category $\C$ to the category of vector spaces over $k$. 
Equivalently, a left $\C$-module $M$ is a collection of $k$-vector spaces $\{ {}_{x}M\}_{x \in \C_{0}}$ provided with a left action 
\[     {}_{y}\C_{x} \otimes {}_{x}M \rightarrow {}_{y}M,     \]
where the image of ${}_{y}f_{x} \otimes {}_{x}m$ is denoted by ${}_{y}f_{x}.{}_{x}m$ or $f m$, satisfying the usual axioms 
\begin{gather*}
   {}_{z}f_{y}.({}_{y}g_{x}.{}_{x}m) = ({}_{z}f_{y}.{}_{y}g_{x}).{}_{x}m,
   \\
   {}_{x}1_{x}.{}_{x}m = {}_{x}m.
\end{gather*}

Right $\C$-modules are defined in an analogous way. 
Also, a \textbf{$\C$-bimodule} is just a $\C^{e}$-module. 

We shall denote ${}_{\C}\mathrm{Mod}$ the category of left $\C$-modules. 
We will always consider left modules unless we say the opposite. 
\end{defi}

The obvious example of $\C$-bimodule is given by the category itself, i.e. ${}_{y}\C_{x}$ for every $x,y \in \C_{0}$. 
We will denote this bimodule by $\C$. 

Also, we notice that given a $\C$-bimodule $M$, we can take the left $\C$-module $\bigoplus_{x \in \C_{0}} {}_{x}M_ {y}$. 
We will denote this left module also by $M$ and similarly for right $\C$-modules. 

Following Mitchell \cite{Mit1}, we recall some facts concerning the category ${}_{\C}\mathrm{Mod}$ of modules over a $k$-linear category $\C$. 
The category ${}_{\C}\mathrm{Mod}$ is abelian, AB5 and cocomplete since the same is true for ${}_{k}\mathrm{Mod}$. 

\begin{defi}
\
\begin{enumerate}
   \item 
   A left $\C$-module $M$ is \textbf{free} if there is a subset $\{ x_{i} : i \in I \}$ of $\C_{0}$ and a natural isomorphism 
   \[     {}_{\place}M \simeq \bigoplus_{i \in I} {}_{\place}\C_{x_{i}}.     \]       

   \item 
   $M$ is \textbf{small} if the functor $\mathrm{Nat}(M , \place)$ preserves coproducts (where $\mathrm{Nat}(M,N)$ denotes the natural transformations 
from $M$ to $N$). 

   \item 
   $M$ is \textbf{finitely generated} if there is an epimorphism $\bigoplus_{i \in I} {}_{\place}\C_{x_{i}} \twoheadrightarrow {}_{\place}M$ for 
   a finite set $I$ and a subset $\{ x_{i} : i \in I \} \subset \C_{0}$. 

   \item 
   $M$ is \textbf {projective} if it is a projective object in the category of $k$-linear functors from $\C$ to ${}_{k}\mathrm{Mod}$.  
\end{enumerate}
\end{defi}

\begin{obs}
\label{proy}
\
\begin{enumerate}

   \item 
   It is easy to see that  a $\C$-module $M$ is projective if and only if it is a retract of a free $\C$-module. 
   In order to obtain this result, let us take $I = \coprod_{x \in \C_{0}} {}_{x}M$. 
   Using Yoneda's lemma:
   \[     \mathrm{Nat} (\bigoplus_{i \in I}{}_{\place}\C_{x_{i}} , M) \simeq \prod_{i \in I} \mathrm{Nat} ({}_{\place}\C_{x_{i}} , M) 
       \simeq \prod_{i \in I} {}_{x_{i}}M.     \]   
   The epimorphism is obtained by considering the collection $(m_{i})_{i \in I}$ with $m_{i} = i$, $\forall i \in I$. 
   Now the remark follows. 

   \item 
   The previous remark assures that the category ${}_{\C}\mathrm{Mod}$ has enough projectives. 

   \item 
   The category ${}_{\C}\mathrm{Mod}$ has also enough injectives. 

   \item
   The category ${}_{\C}\mathrm{Mod}$ has enough generators since this is true in ${}_{k}\mathrm{Mod}$. 
\end{enumerate}
\end{obs}

Next we shall recall the definition of Hochschild-Mitchell homology and cohomology. 
Considering the previous remarks, standard (co)homological methods are available in ${}_{\C}\mathrm{Mod}$. 

\begin{defi}
\label{obstor}
Let $( x_{n+1},\dots,x_{1} )$ be a $(n+1)$-sequence of objects of $\C$. 
The \textbf{$k$-nerve associated to the $(n+1)$-sequence} is the $k$-vector space 
\[     {}_{x_{n+1}}\C_{x_{n}} \otimes \dots \otimes {}_{x_{2}}\C_{x_{1}}.     \]
The \textbf{$k$-nerve of $\C$ in degree $n$} ($n \in \NN_{0}$)  is 
\[     \overline {N}_{n} = \bigoplus_{(n+1)\hbox{\footnotesize-tuples}} {}_{x_{n+1}}\C_{x_{n}} \otimes \dots 
        \otimes {}_{x_{2}}\C_{x_{1}}.     \]      

There is $\C^{e}$-bimodule associated to $\overline{N}_{n}$ defined by  
\[     {}_{y}(N_{n})_{x} = \bigoplus_{(n+1)\hbox{\footnotesize-tuples}} {}_{y}\C_{x_{n+1}} \otimes {}_{x_{n+1}}\C_{x_{n}} 
       \otimes \dots \otimes {}_{x_{2}}\C_{x_{1}} \otimes {}_{x_{1}}\C_{x}.     \]    

Then the associated \textbf{Hochschild-Mitchell complex} is
\[     \dots \stackrel{d_{n+1}}{\rightarrow} N_{n} \stackrel{d_{n}}{\rightarrow} \dots
         \stackrel{d_{2}}{\rightarrow} N_{1} \stackrel{d_{1}}{\rightarrow}
         N_{0} \stackrel{d_{0}}{\rightarrow} \C \rightarrow 0,     \] 
where $d_{n}$ is given by the usual formula, i.e. 
\[     d_{n} (f_{0} \otimes \dots \otimes f_{n+1}) = \sum\limits_{k = 0}^{n} (-1)^{k} f_{0} \otimes \dots \otimes f_{k}.f_{k+1} \otimes \dots \otimes f_{n+1}.     \]

This complex is a projective resolution of the $\C$-bimodule $\C$. 
The proof that it is a resolution is similar to the standard proof for algebras. 
For a detailed proof of the projective part cf. \cite{Hers1}, obs. 4.24, pp. 44-45. 
\end{defi}

\begin{defi}
Given a $\C$-bimodule $M$ the \textbf{Hochschild-Mitchell cohomology of $\C$ with coefficients in $M$} is the cohomology of the following cochain 
complex 
\[     0 \rightarrow \prod_{x \in \C_{0}} {}_{x}M_{x} \stackrel{d^{0}}{\rightarrow} 
        \mathrm{Hom}(N_{1},M) \stackrel{d^{1}}{\rightarrow} \dots \stackrel{d^{n-1}}{\rightarrow} 
       \mathrm{Hom}(N_{n},M) \stackrel{d^{n}}{\rightarrow} \dots,     \]
where $d$ is given by the usual formula, and   
  \begin{align*}
    C^{n}(\C , M) &= \mathrm{Hom}(N_{n} , M) = \mathrm{Nat}(N_{n} , M)      
    \\
    &= \prod_{(n+1)\hbox{\footnotesize-tuples}} \mathrm{Hom}_{k} ({}_{x_{n+1}}\C_{x_{n}} \otimes \dots \otimes 
       {}_{x_{2}}\C_{x_{1}},{}_{x_{n+1}}M_{x_{1}}).    
  \end{align*} 
We denote it $H^{\bullet}(\C , M)$. 

Analogously the \textbf{Hochschild-Mitchell homology of $\C$ with coefficients 
in $M$} is the homology of the chain complex 
\[     \dots \stackrel{d_{n+1}}{\rightarrow} M \otimes N_{n} \stackrel{d_{n}}{\rightarrow} \dots
         \stackrel{d_{2}}{\rightarrow} M \otimes N_{1} \stackrel{d_{1}}{\rightarrow}
       \bigoplus_{x \in \C_{0}} {}_{x}M_{x} \rightarrow 0,     \] 
where $d$ is given by the usual formula and
  \begin{align*}
    C_{n}(\C , M) &= M \otimes_{\C^{e}} N_{n}  
    \\
    &= \bigoplus_{(n+1)\hbox{\footnotesize-tuples}} 
    {}_{x_{1}}M_{x_{n+1}} \otimes {}_{x_{n+1}}\C_{x_{n}} \otimes \dots \otimes {}_{x_{2}}\C_{x_{1}}.     
  \end{align*} 
We denote it $H_{\bullet}(\C , M)$. 
\end{defi}

\begin{obs}
These definitions agree with the usual definition of the Hochschild (co)homology for unital algebras when $\C_{0}$ is finite (cf. \cite{C-R1}, prop. 2.7). 
\end{obs}

One of the main features of Hochschild-Mitchell (co)homology is, as we will next prove, that it is invariant under Morita equivalences. 
It is important to notice the following: 

\begin{teo}
(\cite{C-S1}, thm. 4.7) 
Any Morita equivalence between $k$-linear categories is a composition of equivalences of categories, contractions and expansions. 
\end{teo}

This result leads to the proof of Morita invariance of Hochschild-Mitchell (co)homology. 
However, we shall provide a different proof. 

Observe that two $k$-linear categories are Morita equivalent if and only if their completions (i.e., the additivization of the idempotent completion) are 
Morita equivalent. 
These completions are amenable categories, so they are Morita equivalent if and only if they are equivalent. 
It is then sufficient to prove that Hochschild-Mitchell (co)homology is invariant under: equivalences, additivization and idempotent completion. 
Let us add that the last two processes are compositions of contractions, expansions and equivalences. 
For more details on these facts see the appendix of \cite{C-S1}.  

\begin{defi}
\label{defind}
Given two $k$-linear categories $\C$ and $\D$, a left $\D$-module $N$ and a $k$-linear functor 
$F : \C \rightarrow \D$,
we define the left $\C$-module $FN$ given by
\[     {}_{x}(FN) = {}_{F(x)}N, \forall x \in \C_{0}     \] 
where the action is the following  
\[     f.n = F(f).n,     \]
for $n \in {}_{x}(FN)$ and $f \in {}_{z}\C_{x}$.
The analogous construction is made for right modules and bimodules. 
\end{defi}

\begin{teo}
\label{quasi}
(c.f. \cite{Hers1}, thm. 4.27) 
Let $\C$ and $\D$ be two $k$-linear categories, let $N$ be a $\D$-bimodule and let $F$ be a $k$-linear functor, 
$F : \C \rightarrow \D$,
which is an equivalence of categories. 
Then there are isomorphisms of $k$-modules
  \begin{align*}
    H^{\bullet}(\C , FN) &\simeq H^{\bullet}(\D , N),     
    \\
    H_{\bullet}(\C , FN) &\simeq H_{\bullet}(\D , N).
  \end{align*} 
\end{teo}
\noindent \textbf {Proof.} 
The functor induced by $F$ from ${}_{\D^{e}}\mathrm{Mod}$ to ${}_{\C^{e}}\mathrm{Mod}$ is an equivalence since $F$ is itself an equivalence. 
Hence it is exact and preserves both projective and injective objects. 
Also, $H^{\bullet} (\C , \place)$ is a coeffaceable universal $\delta$-functor. 
As a consequence the following collection of functors 
   \begin{gather*}
     A^{n} : {}_{\D^{e}}\mathrm{Mod} \rightarrow {}_{k}\mathrm{Mod}
     \\ 
     M \mapsto FM \mapsto H^{n}(\C , FM)     
   \end{gather*}
is a universal $\delta$-funtor. 

The collection of functors below is a universal $\delta$-functor too  
  \begin{gather*}
     B^{n} : {}_{\D^{e}}\mathrm{Mod} \rightarrow {}_{k}\mathrm{Mod}
     \\ 
     M \mapsto H^{n}(\D , M).          
   \end{gather*}

In order to prove that two universal $\delta$-functors related by a natural transformations are isomorphic, it is enough to prove that they are isomorphic 
in degree zero (see \cite{Wei1}, section 2.1, p. 32). 
So we are going to prove that $H^{0} (\D , M) \simeq H^{0} (\C , FM)$. 

For this we define the natural morphism 
   \begin{gather*}
     i : H^{0} (\D , M) \rightarrow H^{0} (\C , FM)    
     \\
     ({}_{y}m_{y})_{y \in \D_{0}} \mapsto ({}_{x}im_{x})_{x \in \C_{0}}, 
   \end{gather*}
where ${}_{x}im_{x} = {}_{F(x)}m_{F(x)}$. 

Since $({}_{y}m_{y})_{y \in \D_{0}} \in H^{0}(\D , M)$ we have that ${}_{y'}f_{y}.{}_{y}m_{y} = {}_{y'}m_{y'}.{}_{y'}f_{y}$, where ${}_{y'}f_{y} \in {}_{y'}\D_{y}$. 
Then it also follows that ${}_{x'}f_{x}.{}_{x}im_{x} = {}_{x'}im_{x'}.{}_{x'}f_{x}$, for ${}_{x'}f_{x} \in {}_{x'}\C_{x}$. 

If $y$ is not in the image of $F$, since it is an equivalence, there exists $x \in \C_{0}$ such that $y \simeq F(x)$ with isomorphism $h \in {}_{F(x)}\D_{y}$. 
Taking into account the property mentioned above for the elements of $H^{0}(\D , M)$ we get
\[     {}_{y}m_{y} = h^{-1}.{}_{F(x)}m_{F(x)}.h = h^{-1}.{}_{x}im_{x}.h,     \]
i.e. $i$ is an isomorphism.

The homological case is analogous using the following natural well-defined isomorphism 
   \begin{gather*}
     j : H_{0} (\D , M) \rightarrow H_{0} (\C , FM)    
     \\
     \overline{\sum\limits_{y \in \D_{0}} {}_{y}m_{y}} \mapsto \overline{\sum\limits_{x \in \C_{0}} {}_{x}jm_{x}}, 
   \end{gather*}
where ${}_{x}jm_{x} = {}_{F(x)}m_{F(x)}$. 

The map $j$ is an isomorphism due to the following: if $y$ is not in the image of $F$, since $F$ is an equivalence, 
there exists $x \in \C_{0}$ such that $y \simeq F(x)$ with isomorphism $h \in {}_{F(x)}\D_{y}$. 
Using the definition of $H_{0}(\D , M)$, i.e. $\overline{{}_{y'}f_{y}.{}_{y}m_{y}} = \overline{{}_{y'}m_{y'}.{}_{y'}f_{y}}$, we get 
$\overline{{}_{y}m_{y}} = h^{-1}.\overline{{}_{F(x)}m_{F(x)}}.h = h^{-1}.\overline{{}_{x}im_{x}}.h$. 
\qed

Before stating the following theorem we recall that we have the inclusion functor from $\C$ to its idempotent completion $\hat{\C}$, which we denote 
$hat$, and the inclusion functor from $\C$ to its additivization $\mathrm{Mat}(\C)$, which we denote $mat$. 

\begin{teo}
\label{Morita}
The Hochschild-Mitchell (co)homology is Morita invariant.
\end{teo}
\noindent\textbf{Proof.}
It is enough to prove that Hochschild-Mitchell (co)homology is invariant under additivization and idempotent completion. 

The functor $mat : \C \rightarrow \mathrm{Mat}(\C)$ induces an equivalence of categories 
${}_{\C^{e}}\mathrm{Mod} \simeq {}_{\mathrm{Mat}(\C)^{e}}\mathrm{Mod}$ which we call also $mat$. 
We want to see that 
   \begin{gather*}
     H^{\bullet} (\C , mat M) \simeq H^{\bullet} (\mathrm{Mat}(\C) , M),
   \end{gather*}
for a $\mathrm{Mat}(\C)^{e}$-module $M$, and the same for homology.

Let us define the collections of functors:
\[
\xymatrix@R-20pt
{
  A^{\bullet} : {}_{\mathrm{Mat}(\C)^{e}}\mathrm{Mod} \rightarrow {}_{k}\mathrm{Mod}
  &
  &
  B^{\bullet} : {}_{\mathrm{Mat}(\C)^{e}}\mathrm{Mod} \rightarrow {}_{k}\mathrm{Mod}
  \\
  M \mapsto mat M \mapsto H^{\bullet} (\C , mat M),
  &
  &
  M \mapsto H^{\bullet} (\mathrm{Mat}(\C) , M). 
}
\]

Both are universal $\delta$-functors (coeffaceable in the homological case and effaceable in the cohomological one) since 
the first one is the composition of an exact functor which preserves projectives and injectives (because it is an equivalence) and a universal 
$\delta$-functor; whereas for the second the verification is direct. 

As above we only need to prove that there is a natural isomorphism in degree zero. 
This fact follows as before defining 
  \begin{gather*}
     l : H^{0} (\mathrm{Mat}(\C) , M) \rightarrow H^{0} (\C , mat M)
     \\
     ({}_{(x_{1}, \dots , x_{n})}m_{(x_{1}, \dots , x_{n})})_{(x_{1}, \dots , x_{n}) \in \mathrm{Mat}(\C)_{0}} \mapsto ({}_{x}lm_{x})_{x \in \C_{0}},
   \end{gather*}    
where ${}_{x}lm_{x} = {}_{x}m_{x}$, and
  \begin{align*}
     p : H_{0} (\mathrm{Mat}(\C) , M) &\rightarrow H_{0} (\C , mat M)    
     \\
     \sum\limits_{(x_{1}, \dots , x_{n}) \in \mathrm{Mat}(\C)_{0}} \overline{{}_{(x_{1}, \dots , x_{n})}m_{(x_{1}, \dots , x_{n})}} 
     &\mapsto \sum\limits_{x \in \C_{0}} \overline{{}_{x}pm_{x}}, 
   \end{align*}
where ${}_{x}pm_{x} = {}_{x}m_{x}$. 
It is immediate to see that they are both well-defined and surjective. 

To conclude this proof we remark that the $k$-modules ${}_{(x_{1}, \dots , x_{n})}M_{(x_{1}, \dots , x_{n})}$ over $\mathrm{Mat}(\C)$ can be seen as the 
space of matrices with elements in ${}_{x_{i}}M_{x_{j}}$. 
This is due to the fact that the element $(x_{1}, \dots , x_{n})$ is the coproduct of $x_{1}, \dots, x_{n}$, and then 
$((x_{1} , \dots , x_{n}) , (y_{1} , \dots , y_{m}))$ is the coproduct of $(x_{i}, y_{j})$ 
in $\mathrm{Mat}(\C)^{e}$ ($i = 1, \dots, n$, $j = 1 , \dots , m$), so
\[     M(\bigoplus_{\underset{j = 1 , \dots , m}{i = 1, \dots, n}} (x_{i},y_{j})) \simeq \bigoplus_{\underset{j = 1 , \dots , m}{i = 1, \dots, n}} M((x_{i},y_{j})),     \]
since $M$ es additive.
In our case the maps $p$ y $l$ are injective, since 
\[     e_{i}.{}_{(x_{1}, \dots , x_{n})}m_{(x_{1}, \dots , x_{n})} = {}_{x_{i}}m_{x_{i}}.e_{i},     \]
where $e_{i} \in {}_{x_{i}}\mathrm{Mat}(\C)_{(x_{1}, \dots , x_{n})}$ is defined by $e_{i} = (0 \dots {}_{x_{i}}1_{x_{i}} \dots 0)$. 
Then $l$ and $p$ are isomorphisms. 

For the idempotent completion the proof is even easier. 
The functor $hat : \C \rightarrow \hat{\C}$ induces an equivalence of categories 
${}_{\C^{e}}\mathrm{Mod} \simeq {}_{\hat{\C}^{e}}\mathrm{Mod}$, which we also call $hat$. 
We want to see that 
   \begin{gather*}
     H^{\bullet} (\C , hat M) \simeq H^{\bullet} (\hat{\C} , M),
   \end{gather*}
for any $\hat{\C}^{e}$-module $M$, and the same for homology. 

Let us define the collections of functors:
\[
\xymatrix@R-20pt
{
  A^{\bullet} : {}_{\hat{\C}^{e}}\mathrm{Mod} \rightarrow {}_{k}\mathrm{Mod}
  &
  &
  B^{\bullet} : {}_{\hat{\C}^{e}}\mathrm{Mod} \rightarrow {}_{k}\mathrm{Mod}
  \\
  M \mapsto hat M \mapsto H^{\bullet} (\C , hat M),
  &
  &
  M \mapsto H^{\bullet} (\hat{\C} , M),
}
\]
and analogously for homology. 
Both are universal $\delta$-functors since the first one is the composition of an exact functor which preserves projectives and injectives 
(it is an equivalence) and 
a universal $\delta$-functor; for the second the proof is direct.

Again we only need to prove that there is a natural isomorphism in degree zero. 
We define 
\[
\xymatrix@R-20pt
  {
  l : H^{0} (\hat{\C} , M) \rightarrow H^{0} (\C , hat M)
  &  
  &
  p : H_{0} (\hat{\C} , M) \rightarrow H_{0} (\C , hat M)    
  \\
  ({}_{(x , f)}m_{(x , f)})_{(x , f) \in \hat{\C}_{0}} \mapsto ({}_{x}lm_{x})_{x \in \C_{0}},
  &
  &
  \overline{\sum\limits_{(x , f) \in \hat{\C}_{0}} {}_{(x , f)}m_{(x , f)}} 
  \mapsto \overline{\sum\limits_{x \in \C_{0}} {}_{x}pm_{x}}, 
  }  
\]
where ${}_{x}lm_{x} = {}_{(x, {}_{x}1_{x})}m_{(x , {}_{x}1_{x})}$, and ${}_{x}pm_{x} = {}_{(x, {}_{x}1_{x})}m_{(x , {}_{x}1_{x})}$. 
It is immediate to see that both are well-defined and surjective.  

Moreover, since the objects of $\hat{\C}$ are retracts of objects of $\C$, the elements ${}_{(x , f)}m_{(x , f)}$ are completely determined by the elements of 
the form ${}_{x}m_{x}$ ($x \in \C_{0}$). 
For $f \in {}_{(x , {}_{x}1_{x})}\hat{\C}_{(x , f)}$ and ${}_{(x , f)}m_{(x , f)}$, we get that 
\[     f.{}_{(x , f)}m_{(x , f)} = {}_{x}m_{x}.f.     \]
On the other hand, $f.{}_{(x , f)}m_{(x , f)} = {}_{(x , f)}m_{(x , f)}$, since $f \in {}_{(x , f)}\hat{\C}_{(x , f)}$ and $f = {}_{(x , f)}1_{(x , f)}$. 
As a consequence, $l$ and $p$ are isomorphisms. 

\section{The spectral sequence}

Given a group $G$, an action of $G$ on the $k$-linear category $\C$ is a group morphism from $G$ to $\mathrm{Aut}_{k}(\C) = 
\{ \text{$k$-linear invertible endofunctors of $\C$}  \}$. 
More generally, we define, given a $k$-linear category $\C$ and a $k$-Hopf algebra $(H, \mu, \Delta,\epsilon, \eta)$, 

\begin{defi}
$\C$ is an \textbf{$H$-module category} if each morphism space ${}_{y}\C_{x}$ is an $H$-module, each endomorphism algebra ${}_{x}\C_{x}$ is an 
$H$-module algebra and compositions of maps are morphisms of $H$-modules, where the tensor product of $H$-modules is considered as $H$-module 
via the standard action. 
\end{defi} 

For more details on this definition, see \cite{C-S1}. 

\begin{obs}
Given a  group $G$, it is known that the group algebra $kG$ is a Hopf algebra. 
For a $k$-linear category $\C$, having an action of $G$ is equivalent to the fact of being a $kG$-module category. 
\end{obs}

Let $\C$ be an $H$-module category. 
In \cite{C-S1}, the authors defined the smash product $k$-category as follows: 
the objects of $\C \# H$ are the objects of $\C$, while given two objects $x$ and $y$ in $(\C \# H)_{0} = \C_{0}$, 
${}_{y}(\C \# H)_{x} = {}_{y}\C_{x} \otimes H$. 
Mimicking the smash product of an $H$-module algebra by the Hopf algebra $H$, the composition of maps is given by the following formula: 
\begin{gather*}
   {}_{z}(\C \# H)_{y} \otimes {}_{y}(\C \# H)_{x} \rightarrow  {}_{z}(\C \# H)_{x}
   \\
   ({}_{z}f_{y} \otimes h) \circ ({}_{y}g_{x} \otimes h') =  {}_{z}f_{y} \circ (h_{(1)} {}_{y}g_{x}) \otimes h_{(2)} h',     
\end{gather*}
where we have used Sweedler's notation for the coproduct, i.e. $\Delta(h) = h_{(1)} \otimes h_{(2)}$, $\forall \hskip 1mm h \in H$. 

\begin{obs}
Let $\C$ be a finite $H$-module category. 
Then the $k$-algebras $a(\C) \# H$ and $a(\C \# H)$ are canonically isomorphic. 
\end{obs}

In an analogous way, we may define the notion of comodule categories (see \cite{Mo1} for the definition of $H$-comodule algebra):

\begin{defi}
Given a $k$-linear category $\C$ and a Hopf algebra $H$, we will say that $\C$ is an \textbf{$H$-comodule category} if each morphism space 
${}_{y}\C_{x}$ is an $H$-comodule, each endomorphism algebra ${}_{x}\C_{x}$ is an $H$-comodule algebra and compositions of maps are morphisms of 
$H$-comodules, where the tensor product of $H$-comodules is considered as an $H$-comodule via the standard coaction. 
\end{defi} 

One interesting example is obtained considering a group $G$ and the Hopf algebra $H = kG$ with its usual structure. 
Then our definition of $kG$-comodule category coincides with the definition of $G$-graded $k$-category of \cite{C-M1}. 

\begin{defi}
We say that a $kG$-comodule category $\C$ is \textbf{strongly graded} if 
\[     \sum_{y \in \C_{0}} {}_{z}\C^{s}_{y} \cdot {}_{y}\C^{t}_{x} = {}_{z}\C^{st}_{x}, \forall \hskip 1mm x,z \in \C_{0}, \forall \hskip 1mm s,t \in G,     \]
where we have denoted the composition map by $\cdot$.
\end{defi} 

\begin{ejem}
Let $\C$ be a $k$-linear category provided with an action of a group $G$ and let $\D = \C \# kG$ be the smash product category. 
Then $\D$ is a $kG$-comodule category 
since for all $x, y \in \D_{0} = \C_{0}$, 
\[     {}_{y}\D_{x} = {}_{y}\C_{x} \otimes H = \bigoplus_{g \in G} {}_{y}\C_{x} \otimes k g     \]
and the axioms can be easily verified. 
Moreover, $\D$ is strongly graded. 
\end{ejem}

Next we shall extend to $H$-comodule categories the definition of Galois extension for algebras (see \cite{Mo1}, p. 123). 
Consider a right $H$-comodule category $\D$, a $k$-linear category $\C$ and a $k$-functor $F: \C \rightarrow \D$. 
The category $\D$ is then a left and right $\C$-bimodule. 

The following definition is due to Mitchell (see \cite{Mit2}):
\begin{defi}
Given right and left $\C$-modules $M$ and $N$ respectively, we define their \textbf {tensor product over $\C$}, $M \otimes_{\C} N$, as the $k$-module 
given by
\[     M \otimes_{\C} N = (\bigoplus_{x \in \C_{0}} M_{x} \otimes {}_{x}N) / 
       \cl{\{ m.f \otimes n - m \otimes f.n : m \in M_{x}, n \in {}_{y}N, f \in {}_{x}\C_{y} \}}.     \]
\end{defi} 

\begin{obs}
If $M$ and $N$ are $\C$-bimodules then we can define the $\C$-bimodule \textbf {tensor product over $\C$} by: 
\[     {}_{y}(M \otimes_{\C} N)_{x} = (\bigoplus_{z \in \C_{0}} {}_{y}M_{z} \otimes_{k} {}_{z}N_{x}) / 
       \cl{\{ m.f \otimes n - m \otimes f.n \}}.     \]
for $m \in {}_{y}M_{y'}$, $n \in {}_{x'}N_{x}$, $f \in {}_{y'}\C_{x'}$. 
$M \otimes_{\C} N$ is a $\C$-bimodule in a natural way.
\end{obs}

Given a $kG$-module category $\C$, let us take again $\D = \C \# kG$. 
It is clear that $\D$ is a $\C$-bimodule, by means of the inclusion functor $F : \C \rightarrow \C \# kG = \D$. 
$\D$ is also a $kG$-comodule category. 
We are able to consider $\D \otimes_{\C} \D$, which is a $\C$-bimodule. 

Given $k$-linear categories $\C \subset \D$ such that $\D$ is an $H$-comodule category, we remark that can define the 
$\C$-bimodule $\D \otimes H$ given by ${}_{y}(\D \otimes H)_{x} = {}_{y}\D_{x} \otimes H$ and with action 
$c.(f \otimes h).c' = c.f.c'_{(0)} \otimes h c'_{(1)}$, for $c \in {}_{y'}\C_{y}$, $c' \in {}_{x}\C_{x'}$, $f \in {}_{y}\D_{x}$, $h \in H$.   

\begin{defi}
For $k$-linear categories $\C$ and $\D$ such that $\D$ is an $H$-comodule category and $\C = \D^{coH}$, (i.e. ${}_{y}\C_{x} = {}_{y}\D_{x}^{coH}$, 
$\forall x,y \in \C_{0} = \D_{0}$) 
we shall say that $\C \subset \D$ is an \textbf{$H$-Galois extension} if the natural transformation $\beta : \D \otimes_{\C} \D \rightarrow \D \otimes H$ 
defined by 
\begin{gather*}
   {}_{x}\beta_{y} : {}_{x}(\D \otimes_{\C} \D)_{y} \rightarrow {}_{x}(\D \otimes H)_{y}
   \\
   f \otimes_{\C} g \mapsto (f \circ g_{(0)}) \otimes g_{(1)},
\end{gather*} 
is an isomorphism, for $x,y \in \C_{0}$, $f \in {}_{x}\D_{z}$, $g \in {}_{z}\D_{y}$ and we have used Sweedler's notation for the coaction of $\D$, i.e. 
$\rho (g) = g_{(0)} \otimes g_{(1)}$. 

We easily see that this map is well defined since $\C = \D^{coH}$. 
\end{defi}

We have the following theorem (cf. \cite{U1})
\begin{teo}
\label{graded}
Let $\D$ be a $k$-linear $kG$-comodule category, or equivalently, a $G$-graded category and let $\C = \D^{1}$ be the category of coinvariants 
of $\D$. 
Then $\D$ is strongly graded if and only if $\C \subset \D$ is $kG$-Galois. 
\end{teo}
\noindent\textbf{Proof.}
First note that it is equivalent for $\D$ to be strongly graded, i.e. 
\[     \sum_{y \in \D_{0}} {}_{x}\D^{s}_{y} \cdot {}_{y}\D^{t}_{z} = {}_{x}\D^{st}_{z}, \forall \hskip 1mm s,t \in G     \]
or to satisfy 
\[     \sum_{y \in \D_{0}} {}_{x}\D^{s^{-1}}_{y} \cdot {}_{y}\D^{s}_{x} = {}_{x}\D^{1}_{x}, \forall \hskip 1mm s \in G     \]
since if the latter holds then 
\[     {}_{x}\D^{st}_{z} = {}_{x}\D^{st}_{z} \cdot {}_{z}\D^{1}_{z} = {}_{x}\D^{st}_{z} \cdot (\sum_{y \in \D_{0}} {}_{z}\D^{t^{-1}}_{y} \cdot {}_{y}\D^{t}_{z}) 
                                  = \sum_{y \in \D_{0}} {}_{x}\D^{st}_{z} \cdot {}_{z}\D^{t^{-1}}_{y} \cdot {}_{y}\D^{t}_{z} 
                                  \subset \sum_{y \in \D_{0}} {}_{x}\D^{s}_{y} \cdot {}_{y}\D^{t}_{z} \subset {}_{x}\D^{st}_{z}. 
\]

The map $\beta : \D \otimes_{\C} \D \rightarrow \D \otimes kG$ is given by $f \otimes_{\C} g \mapsto \sum_{s \in G} f \circ g_{s} \otimes s$, where we have denoted 
by $g_{s}$ the $s$ component of $g$. 
To see that ${}_{y}\beta_{x}$ is surjective for all $x,y \in \D_{0}$, it is enough to prove it for ${}_{x}\beta_{x}$ for all $x \in \C$, since $\beta$ is a 
morphism of left $\C$-modules. 

Observe that ${}_{x}\beta_{x}$ is surjective if and only if ${}_{x}1_{x} \otimes s \in \mathrm{Im}(\beta)$, for all $s \in G$. 
This condition is equivalent to the existence of elements $f_{i} \in {}_{x}\D_{y_{i}}$, $g_{i} \in {}_{y_{i}}\D_{x}$ such that 
${}_{x}1_{x} \otimes s = \sum_{t \in G, i} f_{i} (g_{i})_{t} \otimes t$, for all $s \in G$, which is also equivalent to 
$\sum_{i} f_{i} (g_{i})_{s} = {}_{x}1_{x}$ and $\sum_{i} f_{i} (g_{i})_{t} = 0$ for $t \neq s$. 
This last statement is the same as 
\[     \sum_{y \in \D_{0}} {}_{x}\D^{s^{-1}}_{y} \cdot {}_{y}\D^{s}_{x} = {}_{x}\D^{1}_{x}.     \]
Hence, we have that $\D$ being strongly graded is equivalent to $\beta$ being surjective. 

We also need to prove injectivity. 
We will show that if $\D$ is strongly graded, then $\beta$ is injective. 
To prove this we first note that, since 
\[     \sum_{y \in \D_{0}} {}_{x}\D^{s^{-1}}_{y} \cdot {}_{y}\D^{s}_{x} = {}_{x}\D^{1}_{x},     \]
we can write ${}_{x}1_{x} = \sum_{i} c_{s^{-1},i} d_{s,i}$ for $c_{s^{-1},i} \in {}_{x}\D^{s^{-1}}_{y_{i}}$, $d_{s,i} \in {}_{y_{i}}\D^{s}_{x}$ and $s \in G$. 

We define ${}_{y}\alpha_{x} : {}_{y}(\D \otimes kG)_{x} \rightarrow {}_{y}(\D \otimes_{\C} \D)_{x}$ by 
$f \otimes s \mapsto \sum_{i} f c_{s^{-1},i} \otimes_{\C} d_{s,i}$. 
Hence, we have ${}_{y}\alpha_{x} \circ {}_{y}\beta_{x} = \id$, since 
\[     ({}_{y}\alpha_{x} \circ {}_{y}\beta_{x}) (f \otimes_{\C} g) = {}_{y}\alpha_{x} (\sum_{s \in G} f g_{s} \otimes s) 
        = \sum_{i, s \in G} f g_{s} c_{s^{-1},i} \otimes_{\C} d_{s,i} = \sum_{i, s \in G} f  \otimes_{\C} g_{s} c_{s^{-1},i} d_{s,i} = f \otimes_{\C} g
\]
for $f \in {}_{y}\D_{z}$, $g \in {}_{z}\D_{x}$, taking into account that $g_{s} c_{s^{-1},i} \in {}_{z}\D^{1}_{x} = {}_{z}\C_{x}$. 
\qed

If $\C \subset \D$ is an $H$-Galois extension then for each $x \in \D_{0}$ there are $r_{i,x} \in {}_{x}\D_{y_{i}}$ and 
$l_{i,x} \in {}_{y_{i}}\D_{x}$, $i \in I_{x}$ and $I_{x}$ finite, such that if $h,k \in H$, $f \in {}_{y}\D_{x}$, 
$c \in {}_{x}\C_{x}$ we have the following properties (cf. \cite{Ste1}, p. 223):
\begin{enumerate}
\label{properties}
   \item ${}_{x}\beta_{x}(\sum_{i \in I_{x}} r_{i,x}(h) \otimes_{\C} l_{i,x}(h)) = {}_{x}1_{x} \otimes h,$
   
   \item $\sum_{i \in I_{x}} r_{i,x}(hk) \otimes_{\C} l_{i,x}(hk) 
   = \sum_{\underset{j \in I_{y_{i}}}{i \in I_{x}}} r_{i,x}(k) r_{j,y_{i}}(h) \otimes_{\C} l_{j,y_{i}}(h) l_{i,x}(k),$

   \item $\sum_{i \in I_{x}} c r_{i,x}(h) \otimes_{\C} l_{i,x}(h) = \sum_{i \in I_{x}} r_{i,x}(h) \otimes_{\C} l_{i,x}(h) c,$

   \item $\sum_{i \in I_{x}} r_{i,x}(h) \otimes_{\C} l_{i,x}(h)_{(0)} \otimes l_{i,x}(h)_{(1)} 
   = \sum_{i \in I_{x}} r_{i,x}(h_{(1)}) \otimes_{\C} l_{i,x}(h_{(1)}) \otimes l_{i,x}(h_{(2)}),$ 

   \item $\sum_{i \in I_{x}} f_{(0)} r_{i,x}(f_{(1)}) \otimes l_{i,x}(f_{(1)}) = {}_{y}1_{y} \otimes f,$

   \item $\sum_{i \in I_{x}} r_{i,x}(h) l_{i,x}(h) = \epsilon(h),$
 
   \item ${}_{y}\beta'_{x}(f \otimes h) = \sum_{i \in I_{x}} f r_{i,x}(h) \otimes_{\C} l_{i,x}(h)),$
\end{enumerate}
where $\beta'$ is the inverse map of $\beta$. 
The proof of these properties is analogous to the algebraic case. 

For instance to prove the second one, we do the following. 
Taking into account that $\beta$ is injective, the equality holds if and only if we aply $\beta$ to both sides.  
We get 
\[     {}_{x}\beta_{x} (\sum_{i \in I_{x}} r_{i,x}(hk) \otimes_{\C} l_{i,x}(hk)) = {}_{x}1_{x} \otimes hk     \]
and also 
\begin{align*}
   {}_{x}\beta_{x}(\sum_{\underset{j \in I_{y_{i}}}{i \in I_{x}}} r_{i,x}(k) r_{j,y_{i}}(h) \otimes_{\C} l_{j,y_{i}}(h) l_{i,x}(k)) 
   &= \sum_{\underset{j \in I_{y_{i}}}{i \in I_{x}}} r_{i,x}(k) r_{j,y_{i}}(h) (l_{j,y_{i}}(h) l_{i,x}(k))_{(0)} \otimes (l_{j,y_{i}}(h) l_{i,x}(k))_{(1)} 
   \\
   &= \sum_{\underset{j \in I_{y_{i}}}{i \in I_{x}}} r_{i,x}(k) r_{j,y_{i}}(h) l_{j,y_{i}}(h)_{(0)} l_{i,x}(k)_{(0)} \otimes l_{j,y_{i}}(h)_{(1)} l_{i,x}(k)_{(1)} 
   \\
   &= \sum_{i \in I_{x}} r_{i,x}(k) {}_{y_{i}}(1_{y_{i}} \otimes h) l_{i,x}(k) 
   \\
   &= {}_{x}1_{x} \otimes h k 
\end{align*}  
which clearly coincide. 
The other properties are also easy to prove.  

We notice from the isomorphism $\D^{e} \otimes_{\C^{e}} M \simeq \D \otimes_{\C} M \otimes_{\C} \D$ that if $\D$ is flat as left and right $\C$-module, 
then $\D^{e}$ is flat as $\C^{e}$-module (cf. \cite{Ste1}, lemma 2.1). 

Denoting by $R : {}_{\D^{e}}\mathrm{Mod} \rightarrow {}_{\C^{e}}\mathrm{Mod}$ the restriction functor, 
we have also the folowing analogue to lemma 2.2 for the algebraic case considered in \cite{Ste1}:
\begin{lema}
If $\D/\C$ is an extension of $k$-linear categories such that $\D$ is flat as a left and right $\C$-module, then the 
$\delta$-functor $H^{\bullet} (\C , \place) \circ R : {}_{\D^{e}}\mathrm{Mod} \rightarrow {}_{k}\mathrm{Mod}$ is effaceable. 
\end{lema}
\noindent\textbf{Proof.}
We know that the $\delta$-functor $H^{\bullet} (\C , \place)$ is effaceable since it is the same as 
$\mathrm{Ext}^{\bullet}_{\C^{e}} (\C , \place)$. 
Thus we get that $H^{\bullet} (\C , \place) \circ R$ is a $\delta$-functor. 

Let $M$ be an injective $\D$-bimodule. 
Since the functor $R$ is right adjoint to $\D^{e} \otimes_{\C^{e}} \place$, which is exact if and only if $\D^{e}$ is $\C^{e}$-flat, 
we get that $R$ preserves injectives and so $R(M)$ is $\C^{e}$-injective. 
Then the composition $H^{\bullet} (\C , \place) \circ R$ is effaceable. 
\qed

Next proposition is quite similar to prop. 2.3 of Stefan's article 
(cf. \cite{Ste1}), but anyway we present the proof since 
the categorical context may cause some difficulties. 
\begin{prop}
\label{h0}
If $\D/\C$ is $H$-Galois and $\D$ is $\C$-flat as left and right module and $M$ is a $\D^{e}$-module, then there exists a natural structure of right 
$H$-module on $H^{0}(\C , M)$ such that $\mathrm{Hom}_{H} (k , H^{0}(\C , M)) \simeq H^{0}(\D , M)$. 
\end{prop}
\noindent\textbf{Proof.}
We have the following commutative diagram 
\[
\xymatrix 
{
   \mathrm{Hom}_{\D} (\D \otimes_{\C} \D , M) 
   \ar[r]<0.1cm>^{\beta'^{*}}
   \ar^{\nu}[d]
   &
   \mathrm{Hom}_{\D} (\D \otimes H , M) 
   \ar[l]<0.1cm>^{\beta^{*}}
   \ar^{\rho}[d]
   \\
   \mathrm{Hom}_{\C} (\D , M) 
   \ar[r]<0.1cm>^{\omega}
   &
   \mathrm{Hom}_{k} (H , \prod_{x \in \D_{0}} {}_{x}M_{x}) 
   \ar[l]<0.1cm>^{\pi}
}
\]
where the vertical maps are the isomorphisms given by 
\begin{gather*}
   {}_{x}(\nu(f)(d))_{y} = {}_{x}f_{y}(\overline{{}_{x}1_{x} \otimes d}), 
   \\
   {}_{x}(\rho(g)(h))_{x} = {}_{x}g_{x} ({}_{x}1_{x} \otimes h), 
\end{gather*}
for $g \in \mathrm{Hom}_{\D} (\D \otimes H , M)$, $f \in \mathrm{Hom}_{\D} (\D \otimes_{\C} \D , M)$, $d \in {}_{x}\D_{y}$ and $h \in H$. 
The morphisms $\omega$ and $\pi$ are defined in order to make the diagram commutative. 
We have that 
\begin{gather*}
   {}_{x}(\omega(f)(h))_{x} = \sum_{i \in I_{x}} r_{i,x}(h) {}_{y_{i}}f_{x}(l_{i,x}(h)), 
   \\
   {}_{y}(\pi(({}_{x}t_{x})_{x \in \D_{0}})_{x}({}_{y}f_{x}) = ({}_{y}f_{x})_{(0)} {}_{x}t_{x}(({}_{y}f_{x})_{(1)}), 
\end{gather*}

From these formulas one can prove in a direct way that if $f \in \mathrm{Hom}_{\C^{e}} (\D , M)$, then $\omega(f) \in \mathrm{Hom}_{k} (H , M^{\C})$, and 
if $g \in \mathrm{Hom}_{k} (H , M^{\C})$ then $\pi(g) \in \mathrm{Hom}_{\C^{e}} (\D , M)$. 
Thus $\omega$ and $\pi$ induce inverse isomorphims between both subspaces. 

Given $({}_{x}m_{x})_{x \in \C_{0}} \in M^{\C}$ we define the $\C$-bimodule map $\phi_{(m)} : \D \rightarrow M$ by 
${}_{x}f_{y} \mapsto {}_{x}m_{x}.{}_{x}f_{y}$. 
Hence, we get a map $\psi_{(m)} = \omega(\phi_{(m)}) \in \mathrm{Hom}_{k} (H , M^{\C})$, which satisfies 
\[     ({}_{x}f_{y})_{(0)} {}_{y}(\psi_{(m)})_{y}(({}_{x}f_{y})_{(1)}) = {}_{x}m_{x}.{}_{x}f_{y}.     \] 

The $H$-module structure is given by 
\[     {}_{y}(({}_{x}m_{x})_{x \in \D_{0}} \cdot h)_{y} = {}_{y}(\psi_{(m)}(h))_{y} 
       = \sum_{i \in I_{y}} r_{i,y}(h) \hskip 1mm {}_{y_{i}}m_{y_{i}} \hskip 1mm l_{i,y}(h).     \]

We deduce immediately from the second of the properties after thm. \ref{graded} that $(m) \cdot (h k) = ((m)\cdot h) \cdot k$ and hence $H^{0}(\C , M)$ is a 
right $H$-module. 
We will denote the action of $H$ on $H^{0}(\C , M)$ also by $m h$ or $m.h$. 

Also
\[     \mathrm{Hom}_{H} (k , H^{0}(\C , M)) = \mathrm{Hom}_{H} (k , M^{\C}) 
       \simeq \{ ({}_{x}m_{x})_{x \in \D_{0}} \in M^{\C} : m h = \epsilon(h) m, \forall \hskip 1mm h \in H \}     \]
where the last isomorphism is given by the map $f \mapsto f(1)$. 

We need to prove that 
\[     \{ (m) \in M^{\C} : m h = \epsilon(h) m , \forall \hskip 1mm h \in H \} = H^{0}(\D , M).     \]
This is straightforward, since if $(m) \in M^{\C}$ is such that $m h = \epsilon(h) m$, then
\begin{align*}
   {}_{x}m_{x}. {}_{x}f_{y} &= ({}_{x}f_{y})_{(0)} {}_{y}(\psi_{(m)})_{y}(({}_{x}f_{y})_{(1)}) = ({}_{x}f_{y})_{(0)} {}_{y}(({}_{x}m_{x}) \cdot ({}_{x}f_{y})_{(1)})_{y} 
   \\
   &= ({}_{x}f_{y})_{(0)} \epsilon(({}_{x}f_{y})_{(1)}) {}_{y}m_{y} = {}_{x}f_{y} {}_{y}m_{y}.     
\end{align*}
For the other implication, take $(m) \in M^{\D}$ and $h \in H$, then, using the sixth of the properties after thm. \ref{graded}
\[     {}_{x}((m) \cdot h)_{x} = \sum_{i \in I_{x}} r_{i,x}(h) {}_{y_{i}}m_{y_{i}} l_{i,x}(h) 
        = \sum_{i \in I_{x}} r_{i,x}(h) l_{i,x}(h) {}_{x}m_{x} = \epsilon(h) {}_{x}m_{x}.     \]
The naturality is trivial. 
\qed

The proof given by Stefan for algebras (cf. \cite{Ste1}, prop. 2.4) to the following proposition also holds in this case. 
\begin{prop}
If $\D/\C$ is $H$-Galois and $\D$ is $\C$-flat as left and right module and $M$ is a $\D^{e}$-module, then there exists a natural structure of right $H$-module on 
each $H^{n}(\C , M)$ ($n \in \NN_{0}$) which coincides with the one given in the previous proposition. 
\end{prop}

The main theorem of this section is: 
\begin{teo}
Let $\D/\C$ be an $H$-Galois extension such that $\D$ is flat as left and right $\C$-module and let $M$ be a $\D$-bimodule. 
There is a spectral sequence relating Hochschild-Mitchell cohomology of $\C$, $\D$ and Hochschild cohomology of $H$ 
\[     E^{p,q}_{2} = H^{p} (H , H^{q} (\C , M)) \Rightarrow H^{p+q}(\D , M)     \]
which is natural in $M$. 
\end{teo}
\noindent\textbf{Proof.} 
We will apply the Grothendieck spectral sequence. 
In order to do this, let us define  the following three functors: 
\[
\xymatrix@R-20pt
{
  F : {}_{\D^{e}}\mathrm{Mod} \rightarrow {}_{k}\mathrm{Mod}
  &
  &
  F(M) = H^{0} (\D , M),
  \\
  F_{1} : \mathrm{Mod}_{H} \rightarrow {}_{k}\mathrm{Mod}
  &
  &
  F_{1}(N) = \mathrm{Hom}_{H} (k , N),
  \\
  F_{2} : {}_{\D^{e}}\mathrm{Mod} \rightarrow \mathrm{Mod}_{H}
  &
  &
  F_{2}(M) = H^{0} (\C , R(M)),
}
\]
where we have written $R(M)$ for the restriction of $M$ to $\C$-bimodules. 
From now on we will not write it unless necessary. 

Firstly, we get that $F_{1} \circ F_{2} = F$, using proposition \ref{h0}. 

Secondly, given an injective $\D^{e}$-module $M$ we prove that $F_{2}(M) = M^{\C}$ is acyclic for $F_{1}$ as follows: 

There is a $\D^{e}$-monomorphism 
\begin{gather*}
   i : M \rightarrow \prod_{(x_{0},y_{0}) \in \D^{e}_{0}} \mathrm{Hom}_{k} ({}_{(x_{0},y_{0})}\D^{e}_{(\place)} , {}_{x_{0}}M_{y_{0}})     
   \\
   {}_{x_{0}}({}_{y}i_{x}({}_{y}m_{x}))_{y_{0}}({}_{x_{0}}f_{x} \otimes {}_{y}g_{y_{0}}) = {}_{x_{0}}f_{x}. {}_{x}m_{y}. {}_{y}g_{y_{0}},  
\end{gather*}
where the $\D^{e}$-module structure of the right hand side is given by 
${}_{x_{0}}((a \otimes b) {}_{x}\phi_{y})_{y_{0}} (f \otimes g) = {}_{x_{0}}\phi_{y_{0}} (f a \otimes b g)$, for 
${}_{x}\phi_{y} \in \prod_{(x_{0},y_{0}) \in \D^{e}_{0}} \mathrm{Hom}_{k} ({}_{(x_{0},y_{0})}\D^{e}_{(x,y)} , {}_{x_{0}}M_{y_{0}})$, 
$a \in {}_{x'}\D_{x}$, $b \in {}_{y}\D_{y'}$, $f \in {}_{x_{0}}\D_{x'}$ and $g \in {}_{y'}\D_{y_{0}}$. 

Since $M$ is injective, there exists a $\D^{e}$-module $N$ such that 
\[     M \oplus N = \prod_{(x_{0},y_{0}) \in \D^{e}_{0}} \mathrm{Hom}_{k} ({}_{(x_{0},y_{0})}\D^{e}_{(\place)} , {}_{x_{0}}M_{y_{0}}).     \] 
Then it is suficient to prove that $\prod_{(x_{0},y_{0}) \in \D^{e}_{0}} \mathrm{Hom}_{k} ({}_{(x_{0},y_{0})}\D^{e}_{(\place)} , {}_{x_{0}}M_{y_{0}})$ is 
$F_{1}$-acyclic. 
For this we observe that there is an isomorphism 
\[     \mathrm{Hom}_{k}({}_{(x_{0},y_{0})}\D^{e}_{(\place)} , {}_{x_{0}}M_{y_{0}})^{\C} 
       \simeq \mathrm{Hom}_{k}({}_{x_{0}}\D \otimes_{\C} \D_{y_{0}} , {}_{x_{0}}M_{y_{0}})     \]
given by $({}_{x}\phi_{x}) \mapsto \overline{\phi}$, such that $\overline{\phi}(f \otimes_{\C} g) = {}_{x_{0}}(\phi)_{x_{0}} (f \otimes g)$. 
We use this map in order to give $\mathrm{Hom}_{k}({}_{x_{0}}\D \otimes_{\C} \D_{y_{0}} , {}_{x_{0}}M_{y_{0}})$ an $H$-module structure: 
\[     (\psi \cdot h) (f \otimes_{\C} g) = \psi(\sum_{i \in I_{x}} f r_{i,x}(h) \otimes_{\C} l_{i,x}(h) g)     \]
for $f \in {}_{x_{0}}\D_{x}$ and $g \in {}_{x}\D_{y_{0}}$. 
Composing $\alpha$ with $\beta'^{*}$ we get the desired isomorphism. 

We shall prove that $\beta'^{*}$ is $H$-linear, where the $H$-module structure of 
$\mathrm{Hom}_{k}({}_{x_{0}}\D_{y_{0}} \otimes H, {}_{x_{0}}M_{y_{0}})$ is the given by the $H$-structure of $H$: 
\begin{align*}
   (\beta'^{*} (\psi) \cdot h) (f \otimes k) &= (\beta'^{*} (\psi)) (f \otimes h k) = ((\psi \circ \beta') \cdot h) (f \otimes h k)
   \\
   &= \psi (\sum_{i \in I_{y_{0}}} f r_{i,y_{0}}(h k) \otimes_{\C} l_{i,y_{0}}(h k)) 
   = \psi (\sum_{\underset{j \in I_{z_{i}}}{i \in I_{y_{0}}}} f r_{i,y_{0}}(k) r_{j,z_{i}}(h) \otimes_{\C} l_{j,z_{i}}(h) l_{i,y_{0}}(k))
   \\
   &= (\psi \cdot h) (\sum_{i \in I_{y_{0}}} f r_{i,y_{0}}(k) \otimes_{\C} l_{i,y_{0}}(k)) 
   = ((\psi \cdot h) \circ \beta') (f \otimes k) = \beta'^{*} (\psi \cdot h) (f \otimes k),
\end{align*}
for $h, k \in H$, $f \in {}_{x_{0}}\D_{y_{0}}$ and $\psi \in \mathrm{Hom}_{k}({}_{x_{0}}(\D \otimes_{\C} \D)_{y_{0}}, {}_{x_{0}}M_{y_{0}})$. 
Since $\mathrm{Ext}^{\bullet}_{H}(k , \place)$ preserves products, the arguments follows proposition 3.2 of \cite{Ste1}, so we get that 
$F_{2}(M)$ is $F_{1}$-acyclic. 

We obtain that the right derived functors yield a spectral sequence $R^{p} F_{1} (R^{q}F_{2} (M)) \Rightarrow R^{p+q}(M)$. 
But $R^{p}F_{1} (\place) = \mathrm{Ext}_{H}^{p} (k , \place) = H^{p}(H , \place)$ if we consider right $H$-modules as $H$-bimodules using $\epsilon$ to 
define the left action; $R^{q}F_{2} (\place) = H^{q}(\C , \place)$ and $R^{n}F (\place) = H^{n}(\D , \place)$. 
\qed

\begin{coro}
\label{group}
Considering a linear category $\C$ with an action of a group $G$ and $\D = \C \# kG$, the extension $\D/\C$ is $H = kG$-Galois and $\D$ is 
$\C$-flat as left and right module because it is $\C$-free, so for any $\D$-bimodule $M$ we obtain a spectral sequence 
\[     H^{p} (G , H^{q} (\C , M)) \Rightarrow H^{p+q}(\D , M)     \]
which is natural in $M$. 
\end{coro}

\begin{obs}
An analogous result for homology can be obtained with a similar argument for the case $H = kG$. 
\end{obs}

\begin{obs}
It is easy to see that the spectral sequence of \cite{C-R1} is a particular case of the previous result, using that $\C/G$ and $\C \# kG$ are Morita equivalent. 
\end{obs}

\section {Decomposition of the spectral sequence for $H = kG$}

In this section consider a $k$-linear category $\C$ provided of an action of a group $G$. 
The envelopping category $\C^{e}$ is also a $kG$-module category, so we consider $\C^{e} \# kG$.
Let $M$ be a $\C^{e} \# kG$-module, we associate to $M$ the following $\C^{e}$-module 
\begin{gather*}
   M \# kG : \C^{e} \mapsto {}_{k}\mathrm{Mod}
   \\
   {}_{x}(M \# kG)_{y} = \bigoplus_{s \in G} {}_{x}M_{sy}, \forall \hskip 1mm x,y \in \C_{0}     
\end{gather*}
with the induced action of $\C^{e}$ on $M \# kG$.  
We notice that $M \# kG$ is a $\C \# kG$-bimodule. 

Given $g \in G$, let us denote by $\cl{g}$ its conjugacy class, $\cl{G}$ the set of all conjugacy classes in $G$ and $\Z(g)$ the centralizer of $g$ in $G$. 
There are two $\C^{e}$-submodules $Mg$ and $M_{\cl{g}}$ of $M \# kG$ defined respectively by 
\[   {}_{x}(Mg)_{y} = {}_{x}M_{gy}  \hskip 1cm \text{and} \hskip 1cm {}_{x}(M_{\cl{g}})_{y} = \bigoplus_{s \in \cl{g}} {}_{x}M_{sy}, \forall \hskip 1mm x,y \in \C_{0}     \]
with the induced actions of $\C^{e}$. 
 
The proof of the following lemma is straightforward:
\begin {lema} 
\label{acc}
Using the above notations, $M \# kG$ and $M_{\cl{g}}$ are left (and right) $\C^{e} \# kG$-modules with the adjoint action, 
whereas $M.g$ is a $\C^{e} \# k{\cal Z}(g)$-module. 
\end{lema}

Next we shall state a theorem relating the Hochschild-Mitchell cohomology of these modules (for details on the proofs, see \cite{Hers1}, p. 64-80). 
\begin{teo} 
\label{isoHM1}
Given a $k$-linear category  $\C$ which is a $kG$-module category, a $\C^{e} \# kG$-module $M$ and an element $g \in G$ such that $G/\Z(g)$ is finite, 
there is an isomorphism, which is $kG$-linear and natural in $M$:
\[     H^{\bullet}( \C , M_{\cl{g}} )  \simeq  \mathrm{Hom}_{k{\cal Z}(g)} (kG , H^{\bullet}( \C , M.g ).     \] 
\end {teo}
\noindent \textbf{Proof.}
Let us define the collections of functors $F^{n} : {}_{\C^{e} \# kG}\mathrm{Mod} \to {}_{kG}\mathrm{Mod}$, given by the composition of 
\[     M \mapsto M_{\cl{g}} \mapsto H^{n}(\C , M_{\cl{g}}),     \] 
and $E^{n} : {}_{\C^{e} \# kG}\mathrm{Mod} \to {}_{kG}\mathrm{Mod}$, given by the composition  
\[     M \mapsto M.g \mapsto H^{n} (\C , M.g) \mapsto \mathrm{Hom}_{k{\cal Z}(g)} ( kG , H^{n}(\C , M.g) ).     \]      

There is a natural isomorphim in degree $0$ between $F^{\bullet}$ and $E^{\bullet}$ (see lemma \ref{isohm2*} below).  
We shall prove that $F^{\bullet}$ y $E^{\bullet}$ are universal $\delta$-functors. 

For $E^{\bullet}$ notice that $M \mapsto M.g$ is an exact functor. 
It preserves products and injectives, since it sends the injective cogenerator in the category of $\C^{e} \# kG$-modules
(see \cite{Mit1}, pp. 102) $\prod_{x,y \in \C_{0}} \mathrm{Hom}_{k} ({}_{(x,y)}(\C^{e} \# kG)_{(\place)} , C)$ 
into another injective object (where $C$ is an injective cogenerator in the category of $k$-modules). 
Moreover, $M \mapsto H^{\bullet}(\C , M)$ is a $\delta$-functor from ${}_{\C^{e} \# k{\cal Z}(g)}\mathrm{Mod}$ into ${}_{k{\cal Z}(g)}\mathrm{Mod}$, 
and $M \mapsto \mathrm{Hom}_{k{\cal Z}(g)} (kG , M)$ is an exact functor, so the composition is a $\delta$-functor.  
Now, when composing the three functors $E^{\bullet}$ becomes a $\delta$-functor. 

For $F^{\bullet}$, $M \mapsto M_{\cl{g}}$ is exact and it preserves injectives, $H^{\bullet}(\C , \place)$ is a $\delta$-functor from 
${}_{\C^{e} \# kG}\mathrm{Mod}$ to ${}_{kG}\mathrm{Mod}$, then the composition is a $\delta$-functor.

Both are universal $\delta$-functors using that $H^{\bullet}(\C, \place)$ coincides with $\mathrm{Ext}^{\bullet}_{\C^{e}}(\C , \place)$ 
which is a effaceable $\delta$-functor from ${}_{\C^{e} \# kG}\mathrm{Mod}$ to ${}_{kG}\mathrm{Mod}$, for any group $G$ acting on $\C$. 
Then it is universal.

In order to prove that $F^{\bullet}$ is effaceable and then universal we just remark that $(\place)_{\cl{g}}$ is an exact functor, preserving injective 
objects and $H^{\bullet}(\C, \place)$ is an effaceable $\delta$-functor. The argument for $E^{\bullet}$ is similar. 
\qed

\begin{lema}
\label{isohm2*}
Let $\C$ be a $k$-linear category provided of an action of a group $G$, let $M$ be a $\C^{e} \# kG$-module, let $g \in G$ such that 
$G/{\cal Z}(g)$ is finite, and let $\{s_{i} : i = 1, \dots ,n \}$ be a set of representatives of the equivalence classes of 
${\cal Z}(g) \backslash G$. 
There exists a $kG$-linear isomorphism which is natural in $M$
  \begin{align*}
       \phi : (M_{\cl{g}})^{\C} &\stackrel{\sim}{\longrightarrow} \mathrm{Hom}_{k{\cal Z}(g)}(kG , (M.g)^{\C})     
       \\
       ({}_{x}m_{sx}) &\mapsto f,    
  \end{align*}
where $f$ is the ${\cal Z}(g)$-linear extension of $\sum_{i = 1}^{n} k_{i}s_{i} \mapsto k_{i_{0}}s_{i_{0}}.({}_{x}m_{sx})$, 
for $s = s_{i_{0}}^{-1}.g.s_{i_{0}}$.
The $kG$-module structure of $\mathrm{Hom}_{k{\cal Z}(g)}(kG , (M.g)^{\C})$ is the one given by $kG$, i.e. it is the coinduced $kG$-module.
\end{lema}  

Using the previous results we obtain:
\begin{teo}
\label {isoHM2*}
With the same hypotheses of the theorem \ref{isoHM1}, there are natural $kG$-isomorphisms 
\[     H^{\bullet}(\C , M \# G) \simeq \bigoplus_{\cl{g} \in \cl{G}} H^{\bullet}(\C , M_{\cl{g}}) \simeq 
       \bigoplus_{\cl{g} \in \cl{G}} \mathrm{Hom}_{k{\cal Z}(g)} (kG , H^{\bullet}(\C , M.g)).     \]
\end{teo}

\begin{coro}
\label{corHM*}
The spectral sequence of corollary \ref{group} decomposes as a direct sum indexed by the conjugacy classes of $G$, as follows: 
\[     E_{2}^{p,q} = H^{p}(G , H^{q}(\C , M \# G)) \simeq \bigoplus_{\cl{g} \in \cl{G}} H^{p}({\cal Z}(g) , H^{q}(\C , M.g)) \Rightarrow H^{p+q} (\C \# kG, M \# kG).     \]
and the spectral sequence is multiplicative. 
\end{coro}
\noindent \textbf {Proof.} 
It follows using theorem \ref{isoHM2*} and Shapiro's lemma.  
\qed

A similar situation, with less restrictions, holds for homology. 
We shall just state the result since the arguments for the proof are similar to the cohomological case. 
\begin{teo}
\label{isoHM2}
Let $\C$ be a $k$-linear category provided of an action of a group $G$ and let $M$ be a $\C^{e} \# kG$-module. 
Then there is a converging spectral sequence 
\[     E^{2}_{p,q} = H_{p}(G , H_{q}(\C , M \# kG)) \simeq \bigoplus_{\cl{g} \in \cl{G}} H_{p}({\cal Z}(g) , H_{q}(\C , M.g)) 
      \Rightarrow H_{p+q}(\C \# kG, M \# kG).     \]    
\end{teo}

\section {Computations of Hochschild-Mitchell cohomology groups}

The aim of this section is to achieve the computation of some Hochschild-Mitchell cohomology groups of infinite quivers. 
It is well-known that given a quiver $Q$ whose underlying graph is a finite tree, the Hochschild cohomology groups of the algebra $kQ/I$ vanish for 
$n \geq 1$ and $I$ and admissible ideal. 
Let us first give some definitions for linear categories that are analogous to some definitions for associative unitary algebras. 
 
\begin {defi}
Let $\C$ be a $k$-linear category. 
We define the \textbf{(ordinary) quiver associated to $\C$} to be the quiver with set of vertices $Q_{0} = \C_{0}$ and such that the cardinality of 
the set of maps from $x$ to $y$ ($x, y \in \C_{0}$) is equal to $\dim_{k} ({}_{y}\mathrm{rad}(\C)_{x}/{}_{y}\mathrm{rad}^{2}(\C)_{x})$. 
We say the category $\C$ is a \textbf {tree} if the underlying graph of its associated quiver has no cycles.
\end {defi}  

We recall that given a finite dimensional $k$-algebra $A$ and a complete set of orthogonal idempotents of $A$, $\{e_{1}, \dots, e_{n}\}$, the 
$k$-linear category $\C_{A}$ associated to $A$ is such that $(\C_{A})_{0} = \{e_{1}, \dots, e_{n}\}$ and ${}_{e_{j}}(\C_{A})_{e_{i}} = e_{j} A e_{i}$ 
($1 \leq i,j \leq n$). 
Then the algebra associated to $\C_{A}$ is clearly $A$ and $\C_{A}$ is Morita equivalent to the category with only one object and $A$ as set of morphisms 
(cf. \cite{C-S1}). 

\begin{obs}
Notice that the quiver associated to the category $\C_{A}$ is equal to the quiver associated to $A$. 
Hence, we see immediately that the category $\C_{A}$ is a tree if and only if the algebra $A$ is so. 
\end{obs}

Since we are interested in infinite quivers we shall need the following definitions.
\begin{defi}
A \textbf{tower of $k$-modules} is a projective system in the category of $k$-modules, indexed by $\NN$ (with the usual order).
\end{defi}

Given a $k$-linear category $\C$ consider the set $\I$ consisting of the full finite subcategories of $\C$, partially ordered by inclusion. 
\begin{obs}
If $\D$ is a full finite subcategory of $\C$ and $M$ is a $\C$-module, then $M$ is also a $\D$-module using the embedding of $\D$ in $\C$. 
We shall still denote it by $M$. 
The same holds for bimodules. 
\end{obs}

Consider now the set $\A = \{ C^{\bullet} (\D , M), \D \in \I \}$. 
It is clear that it is a projective system with morphisms $\pi^{\bullet}_{\D \leq \E} : C^{\bullet} (\E , M) \rightarrow C^{\bullet} (\D , M)$, such that 
$\pi^{\bullet}_{\E \leq \F} \circ \pi^{\bullet}_{\D \leq \E} = \pi^{\bullet}_{\D \leq \F}$, for all $\D, \E, \F \in \A$ such that $\D \subset \E \subset \F$. 
The maps of complexes $\pi^{\bullet}_{\D} : C^{\bullet}(\D , M) \rightarrow C^{\bullet}$ induced by the inclusion of $\D$ in $\C$ for all $\D$ in $\A$ verify that 
$\pi^{\bullet}_{\D \leq \E} \circ \pi^{\bullet}_{\D} = \pi^{\bullet}_{\E}$. 

\begin{prop}
In the above situation we have 
\begin{enumerate}
   \item  
           $\underset{\leftarrow I}{\lim} C^{\bullet} (\D , M) = C^{\bullet} (\C , M)$,   
   \item  
           $\underset{\leftarrow I}{\lim} C^{\bullet} (\D , \D) = C^{\bullet} (\C , \C)$,

\end{enumerate}
\end{prop}

The proof is straightforward. 

\bigskip

Suppose now that $\C_{0}$ is countable and consider the set $\J$ consisting of a countable collection $\{\D_{i}\}_{i \in \NN}$ of finite full 
subcategories of $\C$ such that $\D_{i} \subset \D_{i+1}, \forall ~ i \in \NN$ and $\cup_{i \in \NN} \D_{i} = \C$ and the projective subsystem of $\A$, 
$\B = \{ C^{\bullet} (\D , M), \D \in \J \}$. 
Then $\B$ is a tower of $k$-modules and it is cofinal in $\A$. 
As a consequence 
\[     \underleftarrow {\lim}_{\J} C^{\bullet} (\D , M) = \underleftarrow {\lim}_{\I} C^{\bullet} (\D , M) = C^{\bullet} (\C , M).     \]
Also 
\[     \underleftarrow {\lim}_{\J} C^{\bullet} (\D , \D) = \underleftarrow {\lim}_{\I} C^{\bullet} (\D , \D) = C^{\bullet} (\C , \C).     \]

\begin {defi}
Given a poset $(A , \leq)$ and a projective system indexed by $A$, $(C_{a} , f_{a \leq b})$, we shall say that $(C_{a} , f_{a \leq b})$ satisfies the 
\textbf {Mittag-Leffler condition} (M-L condition) if: 
\[     \forall ~ a \in A, \exists ~ b \geq a \hskip 1mm \text{such that} \hskip 1mm \mathrm{Im}(f_{a \leq b} : C_{b} \rightarrow C_{a}) = \mathrm{Im}(f_{a \leq d} : C_{d} \rightarrow C_{a}), 
\forall ~ d \geq b.     \] 
We shall say  that $(C_{a} , f_{a \leq b})$ satisfies the \textbf {trivial Mittag-Leffler condition} if: 
\[     \forall ~ a \in A, \exists ~ b \geq a \hskip 1mm \text{such that} \hskip 1mm \mathrm{Im}(f_{a \leq d} : C_{d} \rightarrow C_{a}) = 0, \forall ~ d \geq b.     \] 
\end {defi}

Next we recall a result from \cite{Wei1} (prop. 3.5.7, p.83): 
\begin {prop}
Given a tower of $k$-modules $(A_{i} , f_{i \leq j})$ which satisfies the M-L condition, we have that $\underleftarrow {\lim}^{1}_{\NN} A_{i} = 0$, 
where $\underleftarrow {\lim}^{1}_{\NN}$ denotes the first derived functor of $\underleftarrow {\lim}_{\NN}$. 
\end {prop}

The proof of the homological analogue of the following theorem can be found in \cite{Wei1} (teo. 3.5.8, pp. 83-84). 
The cohomological case can be proved in the same way: 
\begin {teo}
\label{lim}
Consider a tower of cochain complexes $(C^{\bullet}_{i} , f^{\bullet}_{i \leq j})$ which satisfies the M-L condition. 
Let $C^{\bullet} = \underleftarrow {\lim}_{\NN} (C^{\bullet}_{i},f^{\bullet}_{i})$. 
There is a short exact sequence 
\[     0 \rightarrow \underleftarrow {\lim}^{1}_{\NN} H^{\bullet -1} (C_{i}) \rightarrow H^{\bullet}(C) 
          \rightarrow \underleftarrow {\lim}_{\NN} H^{\bullet}(C_{i}) \rightarrow 0.     \]
\end {teo}

\begin{obs}
The fact that the category is countable is necessary. 
If it is not the case, $\underleftarrow {\lim}^{n}_{I}$ is no more trivial for $n \geq 2$. 
Suppose for example that the cardinality of $\C$ is $\aleph_{n}$ ($n \geq 0$) then $\underleftarrow {\lim}^{j}_{\I} = 0$ if $j \geq n + 2$ 
(see \cite{Mit2}, section 16, pp. 68-70).
\end{obs}

The above theorem can be applied to Hochschild-Mitchell  cohomology, obtaining in this situation 

\begin {coro}
\label{cohomcat}
Given a $k$-linear countable category $\C$ and a $\C$-bimodule $M$, there are short exact sequences  
\[     0 \rightarrow \underleftarrow {\lim}^{1}_{\J} H^{\bullet -1}(\D , M) \rightarrow H^{\bullet}(\C , M) 
          \rightarrow \underleftarrow {\lim}_{\J} H^{\bullet}(\D , M) \rightarrow 0.     \]
and 
\[     0 \rightarrow \underleftarrow {\lim}^{1}_{\J} HH^{\bullet -1}(\D) \rightarrow HH^{\bullet}(\C) 
         \rightarrow \underleftarrow {\lim}_{\J} HH^{\bullet}(\D) \rightarrow 0.     \]
\end {coro}

Now we can prove the following result. 

\begin {teo}
\label{coarbol}
If $\C$ is a $k$-linear countable category of locally finite dimension, such that $\C$ is a tree, then $HH^{n}(\C) = 0$, 
for all $n \geq 1$.
\end {teo}
\noindent\textbf {Proof.}
If a category is a tree then is is easy to see that every subcategory is also a tree. 
Since for a finite subcategory $\D$, being a tree, we know that $HH^{n}(\D) =0$, for $n \geq 1$, we obtain from thm. \ref{cohomcat} that 
\[     \underleftarrow {\lim}^{1}_{\J} HH^{n - 1}(\D) = HH^{n}(\C), \hskip 1mm \text{for} \hskip 1mm n \geq 1.     \]
As a consequence, $HH^{n}(\C)$ is trivial for $n \geq 2$. 
When $n = 1$, $HH^{0}(\D) = {\cal Z} (\D)$, the center of the linear subcategory $\D$, then   
\[     \underleftarrow {\lim}^{1}_{\J} {\cal Z} (\D) \simeq HH^{1}(\C).     \]
Since $\D_{0}$ is finite and $\C$ is of locally finite dimension, the algebra associated to $\D$ is finite dimensional. 
So ${\cal Z} (\D)$ is also a finite dimensional algebra. 
Then the system $\{{\cal Z} (\D)\}_{\D \in \J}$ satisfies the Mittag-Leffler condition (considering the dimensions), which implies that  
$\underleftarrow {\lim}^{1}_{\J} {\cal Z} (\D) = 0$, and so $HH^{1}(\C) = 0$. 
\qed

Next we shall consider other examples of infinite quivers associated to partially ordered sets. 

\begin{ejem}
\begin{enumerate}
   \item
   Consider the quivers associated to the poset $P^{\infty}_{n} = \NN_{0} \times \{0, \dots, n-1\}$ 
\[
\xymatrix@C-20pt 
   {
   \vdots
   &
   &
   \vdots
   &
   &
   \vdots
   &
   &
   \dots
   &
   \vdots
   &
   &
   \vdots
   \\
   \bullet 
     \ar[d]
     \ar[rrrrrrrrrd]
   &
   (2,0)
   &
   \bullet 
      \ar[d]
      \ar[lld]
   &
   (2,1)
   &
   \bullet 
      \ar[d] 
      \ar[lld] 
   &
   (2,2)
   &
   \dots
   &
   \bullet 
      \ar[d]
      \ar[lld]
   & 
   (2,n-2)
   &
   \bullet 
      \ar[d]
      \ar[lld]
   &
   (2,n-1)
   \\
   \bullet 
     \ar[d]
     \ar[rrrrrrrrrd]
   &
   (1,0)
   &
   \bullet 
      \ar[d]
      \ar[lld]
   &
   (1,1)
   &
   \bullet 
      \ar[d] 
      \ar[lld] 
   &
   (1,2)
   &
   \dots
   &
   \bullet 
      \ar[d]
      \ar[lld]
   & 
   (1,n-2)
   &
   \bullet 
      \ar[d]
      \ar[lld]
   &
   (1,n-1)
   \\
   \bullet 
   &
   (0,0)
   &
   \bullet 
   &
   (0,1)
   &
   \bullet
   &
   (0,2) 
   &
   \dots
   &
   \bullet
   &
   (0,n-2)
   &
   \bullet
   & 
   (0,n-1) 
   }
\]   
and $P^{\infty,\infty}_{n} = \ZZ \times \{0, \dots, n-1\}$
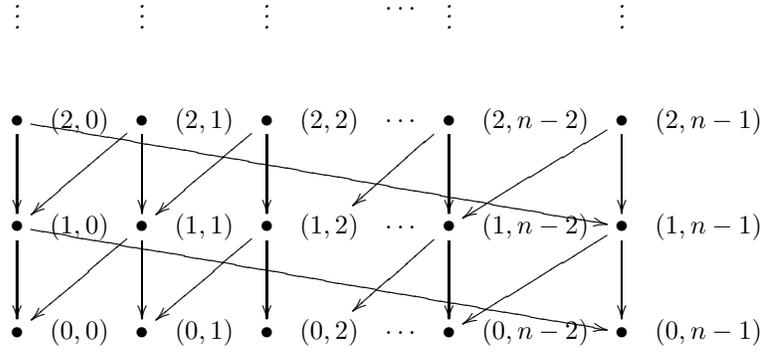
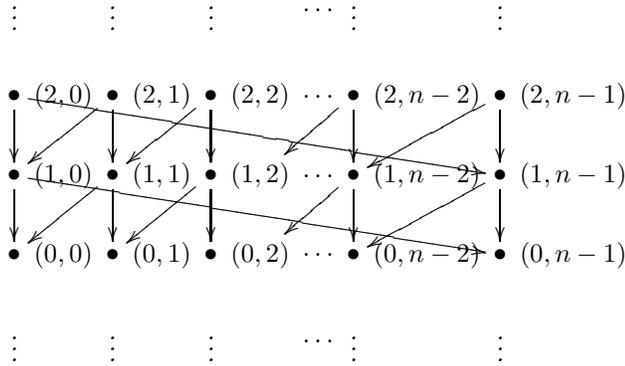
\[
\xymatrix@C-25pt @R-10pt 
   {
   \vdots
   &
   &
   \vdots
   &
   &
   \vdots
   &
   &
   \dots
   &
   \vdots
   &
   &
   \vdots
   \\
   \bullet 
     \ar[d]
     \ar[rrrrrrrrrd]
   &
   (2,0)
   &
   \bullet 
      \ar[d]
      \ar[lld]
   &
   (2,1)
   &
   \bullet 
      \ar[d] 
      \ar[lld] 
   &
   (2,2)
   &
   \dots
   &
   \bullet 
      \ar[d]
      \ar[lld]
   & 
   (2,n-2)
   &
   \bullet 
      \ar[d]
      \ar[lld]
   &
   (2,n-1)
   \\
   \bullet 
     \ar[d]
     \ar[rrrrrrrrrd]
   &
   (1,0)
   &
   \bullet 
      \ar[d]
      \ar[lld]
   &
   (1,1)
   &
   \bullet 
      \ar[d] 
      \ar[lld] 
   &
   (1,2)
   &
   \dots
   &
   \bullet 
      \ar[d]
      \ar[lld]
   & 
   (1,n-2)
   &
   \bullet 
      \ar[d]
      \ar[lld]
   &
   (1,n-1)
   \\
   \bullet 
   &
   (0,0)
   &
   \bullet 
   &
   (0,1)
   &
   \bullet
   &
   (0,2) 
   &
   \dots
   &
   \bullet
   &
   (0,n-2)
   &
   \bullet
   & 
   (0,n-1) 
   \\
   \vdots
   &
   &
   \vdots
   &
   &
   \vdots
   &
   &
   \dots
   &
   \vdots
   &
   &
   \vdots
   }
\] 
   (where $(l,t) < (l+1,t)$, $(l,t) < (l+1,t+1)$ and $(l,n) = (l,0)$)  
   and two $k$-linear categories, $\C_{n}^{\infty}$ and $\C_{n}^{\infty, \infty}$ such that the quivers associated to these categories are respectively 
   $P^{\infty}_{n}$ and $P^{\infty, \infty}_{n}$. 
   Both categories are filtered by full finite subcategories corresponding to algebras of type $A^{0}_{q n + s}$. 
   Hochschild cohomology of this algebra was computed in \cite{G1}, thm. 2.4. 
   The authors obtain the following results ($n \geq 3$, $q \geq 0$, $0 \leq s < n$, $qn+s > 0$): 
   \[     HH^{i}(A^{0}_{q n + s}) = \begin{cases} 
                                                                    k & \text{if $i=0$}
                                                                    \\
                                                                    k & \text{if $i = 2 q + 1$}
                                                                    \\
                                                                    0 & \text{otherwise}
                                               \end{cases}   
   \]   
   for $s \neq 0$, and    
   \[    HH^{i}(A^{0}_{q n}) = \begin{cases}
                                                             k & \text{if $i=0$}
                                                             \\
                                            k^{n-1} & \text{if $i = 2 q$}
                                            \\
                                            0 & \text{otherwise}
                                      \end{cases}                           
   \]       
   Looking at the short exact sequence of theorem \ref{cohomcat}:    
   \[     0 \rightarrow \underleftarrow {lim}^{1}_{\J} HH^{\bullet -1}(\D) \rightarrow HH^{\bullet}(\C^{\infty}_{n})          
             \rightarrow \underleftarrow {lim}_{\J} HH^{\bullet}(\D) \rightarrow 0.     \]   
   We remark that there exists $n_{0}$ such that $HH^{m}(\C^{\infty}_{n}) = 0$ for $m \geq n_{0}$.  
   Hence $\underleftarrow {lim}_{\J} HH^{m}(\D) = 0$ and $\underleftarrow {lim}^{1}_{\J} HH^{m-1}(\D) = 0$, so $HH^{m}(\C^{\infty}_{n}) = 0$ for $m \geq 2$. 

   For $m = 1$ we have that the left term of the short exact sequence is $\underleftarrow {lim}^{1}_{\J} HH^{0}(\D)$. 
   It is zero since $HH^{0}(\D) = k$ and the Mittag-Leffler condition holds.  
   Then $HH^{1}(\C^{\infty}_{n}) = \underleftarrow {lim}^{1}_{\J} HH^{0}(\D) = 0$. 
   Also $HH^{0}(\C^{\infty}_{n}) = \underleftarrow {lim}_{\J} HH^{0}(\D) = k$. 
   The same reasons apply to obtain that $HH^{m}(\C^{\infty,\infty}_{n}) = 0$ for $m \geq 1$ and $HH^{0}(\C^{\infty,\infty}_{n}) = k$.    
   
   \begin{obs}
   Since $\ZZ$ acts on this category we are able to apply our spectral sequence in order to compute $HH^{\bullet}(\C^{\infty, \infty}_{n} \# k\ZZ)$. 
   \end{obs}

   \item
   Consider the category $\C_{n}^{\infty > }$ associated to the poset 
   $P^{\infty >}_{n} = \NN_{0} \times \{0, \dots, n-1\} \cup \{(0,p) : -t \leq p < 0 \} \cup \{ (0,p) : n-1 < p \leq n-1+t'\} \cup \{ (1,p) : -t \leq p \leq n-1+t'\}$ 

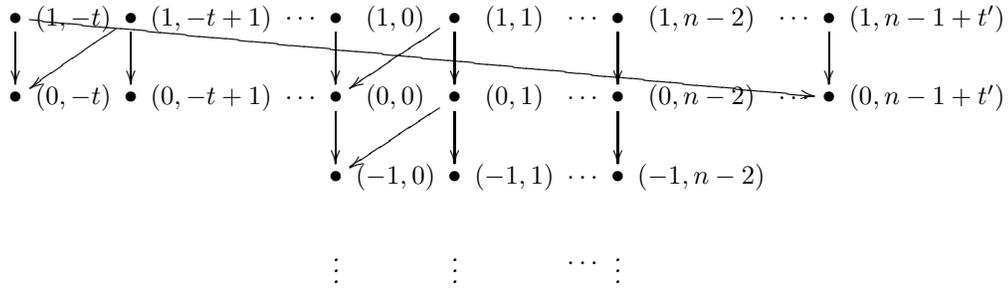
\[
\xymatrix@C-25pt @R-10pt 
   {
   \bullet 
     \ar[d]
     \ar[rrrrrrrrrrrrrd]
   &
   (1,-t)
   &
   \bullet 
      \ar[d]
      \ar[lld]
   &
   (1,-t+1)
   &
   \dots
   &
   \bullet 
      \ar[d] 
   &
   (1,0)
   &
   \bullet 
      \ar[d]
      \ar[lld]
   &
   (1,1)
   &
   \dots
   &
   \bullet 
      \ar[d]
   & 
   (1,n-2)
   &
   \dots
   &
   \bullet 
      \ar[d]
   &
   (1,n-1+t')
   \\
   \bullet 
   &
   (0,-t)
   &
   \bullet 
   &
   (0,-t+1)
   &
   \dots
   &
   \bullet 
      \ar[d] 
   &
   (0,0)
   &
   \bullet 
      \ar[d]
      \ar[lld]
   &
   (0,1)
   &
   \dots
   &
   \bullet 
      \ar[d]
   & 
   (0,n-2)
   &
   \dots
   &
   \bullet 
   &
   (0,n-1+t')
   \\
   &
   &
   &
   &
   &
   \bullet 
   &
   (-1,0)
   &
   \bullet 
   &
   (-1,1)
   &
   \dots
   &
   \bullet 
   & 
   (-1,n-2)
   &
   &
   &
   &
   &
   \\
   &
   &
   &
   &
   & 
   \vdots 
   &
   &
   \vdots 
   &
   &
   \dots
   &
   \vdots
   &
   &
   &
   &
   &  
   }
\]

   (where $(l,l') < (l+1,l')$, $(l,l') < (l+1,l'+1)$, $(l,n) = (l,0)$ for $l \leq 0$ and (1, -t) = (1, n+ t')).   
   It is filtered by full finite subcategories corresponding to algebras of type $A^{>}_{q n}$, whose  
   Hochschild cohomology groups are computed in \cite{G-R1}, thm. 2.1. 
   They have the following results ($n \geq 3$, $q \geq 0$, $0 \leq s < n$, $qn+s > 0$): 
   \[     HH^{i}(A^{>}_{q n}) = \begin{cases} 
                                                                    k & \text{if $i=0$},
                                                                    \\
                                                                    k & \text{if $i = 1$},
                                                                    \\
                                                                    0 & \text{otherwise}.
                                               \end{cases}   
   \]   
   The short exact sequence of our theorem shows that $HH^{m}(\C^{\infty >}_{n}) = 0$ for $m \geq 2$ and $HH^{1}(\C^{\infty >}_{n}) = k$. 
   In order to compute $HH^{0}(\C^{\infty >}_{n})$, it is sufficient to know that the maps of the projective system $(HH^{0}(\D))_{\D}$ are surjective, then 
   $HH^{m}(\C^{\infty >}_{n}) = \underleftarrow{\lim} HH^{0}(\D) = k$.   

   \item
   Consider the category $\C(U)_{\infty}^{n}$ associated to the poset 
   ${}_{\infty}^{n}U = \NN_{0} \times \{0, \dots, n - 1\}$ 

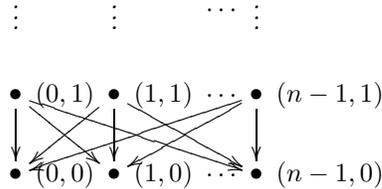
   \[
   \xymatrix@C-25pt @R-10pt 
   {
   \vdots
   &
   &
   \vdots
   &
   &
   \dots
   &
   \vdots
   &
   \\
   \bullet 
     \ar[d]
     \ar[rrrrrd]
     \ar[rrd]
   &
   (0,1)
   &
   \bullet 
      \ar[d]
      \ar[lld]
      \ar[rrrd]
   &
   (1,1)
   &
   \dots
   &
   \bullet 
      \ar[d]
      \ar[llld]
      \ar[llllld]
   & 
   (n-1,1)
\\
   \bullet 
   &
   (0,0)
   &
   \bullet 
   &
   (1,0)
   &
   \dots
   &
   \bullet 
   & 
   (n-1,0) 
   }
\] 
   (where $(l,t) < (l',t+1)$).  
   Again in this case $\C(U)_{\infty}^{n}$ is filtered by full finite subcategories corresponding to algebras of type ${}^{n}_{m}U$, whose  
   Hochschild cohomology groups are computed in \cite{G-R1}, see thm. 2.8. 
   The authors obtain the following results: 
   \[     HH^{i}({}^{n}_{m}U) = \begin{cases} 
                                                                    k & \text{if $i=0$},
                                                                    \\
                                                                    (k^{n-1})^{m + 1} & \text{if $i = m$},
                                                                    \\
                                                                    0 & \text{otherwise}.
                                               \end{cases}   
   \]   
   The short exact sequence of our theorem shows that $HH^{m}(\C(U)_{\infty}^{n}) = 0$ for $m \geq 1$ and $HH^{0}(\C(U)_{\infty}^{n}) = k$. 
\end{enumerate}

\end{ejem}

\begin{obs}
The following fact is well known in representation theory 
(see for example \cite{G-R1}, section 3). 
Let $A$ be an algebra presented as $(Q_{A}, I)$. 
If we have a source $x$ in the quiver $Q_{A}$ then the subcategory 
$A_{x}$ consisting of the objects of $Q_{A}$ except $x$ has 
$Q_{A_{x}}$ as quiver, 
obtained by deleting $x$ and all arrows starting at $x$. 
The presentation of the algebra $A$ yields an induced presentation 
$(Q_{A_{x}}, I')$ of $A_{x}$. 
If $P_{x}$ denotes the projective $A$-module in $x$ and we define 
$M = \mathrm{rad}(P_{x})$, then $A$ is isomorphic to the one point extension 
algebra $A_{x}[M]$. 
Dually, if $x$ is a sink, we take $M = I_{x}/\mathrm{soc}(I_{x})$ for $I_{x}$ 
the injective $A$-module in $x$, and then $A = [M]A_{x}$. 
From Happel's long exact sequence we can see that if all algebras are 
triangular, then $HH^{\bullet}(A_{x}) = HH^{\bullet}(A)$. 

Now we are able to  generalize this method for an infinite sequence of sinks 
and sources:
we suppose that we have a category $\C$ such that its quiver can be decomposed 
into a finite quiver $Q_{A}$ with relations corresponding to an 
algebra $A$ and an infinite sequence of sinks or sources connected to it. 
By taking the obvious countable family of finite subquivers, 
we  see immediately that $HH^{\bullet}(\C) = HH^{\bullet}(A)$. 
\end{obs}



\footnotesize \noindent E.H.:
\\Departamento de Matem\'atica,
 Facultad de Ciencias Exactas y Naturales,\\
 Universidad de Buenos Aires
\\Ciudad Universitaria, Pabell\'on 1\\
1428, Buenos Aires, Argentina. \\
 
\noindent A.S.:
\\Departamento de Matem\'atica,
 Facultad de Ciencias Exactas y Naturales,\\
 Universidad de Buenos Aires
\\Ciudad Universitaria, Pabell\'on 1\\
1428, Buenos Aires, Argentina. \\{\tt asolotar@dm.uba.ar}


\begin{thebibliography}{99}

\bibitem [B-G]{B-G1} Bongartz, D. J.; Gabriel, P. \textit {Covering spaces in representation theory}. Invent. Math. {\bf 65}, 3, (1981-1982), pp. 331--378. 

\bibitem [C-M]{C-M1} Cibils, C.; Marcos, E. \textit {Skew category, Galois covering and smash product of a category over a ring}. To appear in Proc. of the Am. Math. Soc. 
\texttt {http://arxiv.org/abs/math.RA/0312214}.

\bibitem [C-R]{C-R1} Cibils, C.; Redondo M. J. \textit {Cartan-Leray spectral 
sequence for Galois coverings of categories}. 
J. of Alg. {\bf 284}, (2005), pp. 310–-325. 

\bibitem [C-S]{C-S1} Cibils, C.; Solotar, A. \textit {Galois coverings, Morita equivalence and smash extensions of categories over a field}. 
\texttt {http://arxiv.org/abs/math.RA/0502308}.

\bibitem [D-S]{D-S1} Dowbor, P.; Skowronski, A. \textit {Galois coverings of 
representation-infinite algebras}. Comm. Math. Helv. {\bf 62}, (1987), 
pp. 311--337. 


\bibitem [Ga]{Ga1} Gabriel, P. \textit {The universal cover of a representation-finite algebra}. Lect. Notes in Math. {\bf 903}, (1981), Springer-Verlag, pp. 68--105.

\bibitem [G]{G1} Gatica, A. \textit {Cohomolog\'{\i}a de Hochschild de \'algebras de incidencia}. Tesis de Doctorado, Universidad Nacional del Sur, 2001. 

\bibitem [G-R]{G-R1} Gatica, A.; Rey, A. \textit {Hochschild cohomology groups of incidence algebras associated to reduced posets and their
fundamental groups}. To appear in Comm. in Alg. 

\bibitem [Gro]{Gro1} Grothendieck, A. \textit {Sur quelques points 
d'alg`ebre homologique}. Tohoku Math. J. {\bf 9}, (1957), pp.119--221.    

\bibitem [H]{Hers1} Herscovich Ramoneda, E. \textit{La homolog\'{\i}a de 
Hochschild-Mitchell de categor\'{\i}as lineales y su parecido a la 
homolog\'{\i}a de 
Hochschild de \'algebras}. 
Tesis de Licenciatura, Universidad de Buenos Aires, 2005. 

\bibitem [K]{Kel1} Keller, B. \textit {On the cyclic homology of exact 
categories}. J. of Pure and App. Alg. {\bf 136}, 1, (1999), pp. 1--56.

\bibitem [McC]{Mc1} McCarthy, R. \textit {The cyclic homology of an exact 
category}. J. of Pure and App. Alg. {\bf 93}, 3, (1994), pp. 251--296.
 
\bibitem [M1]{Mit1} Mitchell, B. \textit {Theory of categories}. Academic Press Inc., 1965.

\bibitem [M2]{Mit2} Mitchell, B. \textit {Rings with several objects}. 
Adv. in Math. {\bf 8}, (1972), pp.1--161.

\bibitem [Mo]{Mo1} Montgomery, S. \textit {Hopf algebras and their actions on rings}. CBMS, 82, AMS, 1993. 

\bibitem [P-R]{P-R1} Pirashvili, T.; Redondo, M. 
\textit {Abelian categories}. \texttt {http://arxiv.org/abs/math.CT/0504282}.

\bibitem [R]{Re1} Redondo, M.J. \textit {Hochshcild cohomology: some methods for computations}. Resenhas IME-UMP {\bf 5(2)}, (2001), pp. 113--137.

\bibitem [S]{Ste1} Stefan, D. \textit {Hochschild cohomology on Hopf Galois 
extensions}. J. of Pure and App. Alg. {\bf 103}, (1995), pp. 221--233.

\bibitem [Ta]{Ta1} Takeuchi, M. \textit {Morita theorems for categories of 
comodules}. J. Fac. Sci. Univ. Tokyo {\bf 24}, (1977), pp. 1483--1528.

\bibitem [U]{U1} \"Ulbrich, K. \textit {Vollgraduierte Algebren}. 
Abh. Math. Sem. Univ. Hamburg {\bf 51}, (1981), pp. 136--148.

\bibitem [W]{Wei1} Weibel, C. \textit {An introduction to homological algebra}. Cambridge Studies in Advanced Mathematics, 38. 
Cambridge University Press, Cambridge, 1995.  

\bibitem [Wi]{Wi1} Witherspoon, S. \textit{Products in Hochschild cohomology 
and Grothendieck rings of group crossed products}. Adv. in Math. {\bf 185}, 
1, (2004), pp. 136--158. 

\end{thebibliography}
\end{document}